# Exploring the potential resource integration under passenger-freight shared mobility: collaborative optimization of multi-type bus scheduling and dynamic vehicle capacity allocation for urban-rural bus routes

Jiabin Wu[a], Zijian Huang[b], Linhong Wang[c], Yiming Bie[c,]*, Yuting Ji[c], Jun Gong[d]

[a] *School of Management, Foshan University, Foshan 528000, China*
[b] *School of Artificial Intelligence, South China Normal University, Guangzhou 510631, China*
[c] *School of Transportation, Jilin University, Changchun 130022, China*
[d] *Department of Civil Engineering, The University of Hong Kong, Hong Kong 999077, China*

*Corresponding author: **Yiming Bie** (E-mail address: yimingbie@126.com).

**Abstracts**：Under the global background of developing urban-rural travel patterns, traditional urban-rural public transport systems are generally faced with the serious challenges of passenger loss and operating deficit, leading to a reduction in the bus frequency and service reliability. In order to break the vicious circle of 'demand decline-supply shrinkage', passenger-freight shared mobility (PFSM), an innovative operation mode, can achieve synergies between urban-rural logistics and public transport services by integrating public transit network resources and vehicle spare capacity. However, PFSM has changed the operating characteristics of urban-rural bus systems, posing some new challenges. To expand the relevant theory and find the solutions to those challenges, this study proposes an 'economy-efficiency-low-carbon' -oriented resource reconfiguration strategy by formulating the collaborative bilevel optimization of multi-type bus scheduling and dynamic vehicle capacity allocation for urban-rural bus routes. The improved jellyfish search algorithm is developed to solve the premature convergence problem of the traditional algorithms in solving a high-dimensional hybrid discrete-continuous optimization. The results of a case of two urban-rural bus lines in Shanxi Province, China, indicate that the proposed scheme can improve operating revenue by 328.45% and reduce freight carbon emissions by 19.12 tons/year within the increase of 19.46% in average passenger travel time. The sensitivity analysis explicates key parameters selected for PFSM in terms of economic, efficiency and environmental dimensions. The proposed method provides some novel insights and solutions for the sustainable development of urban-rural public transport systems and the last kilometer problem of rural logistics, with significant values of both economic growth and environmental carbon reduction.

# 1. Introduction

Urban-rural public transport is an important link connecting economic exchanges between counties, towns, and rural areas, which is of great significance in promoting rural-urban migration and achieving rural revitalization [1][2]. However, due to the low population density, scattered passenger flow, and long travel distance in rural areas, there exist some common problems on urban-rural bus routes such as long travel time and high fares, which greatly reduces the motivation of residents to travel by bus. In order to avoid serious losses, some public transport enterprises have to reduce bus frequencies, raise fares, and even cancel some routes, leading to the vicious circle of 'demand shrinking-service degradation' [3]. Moreover, the rapid development of rural e-commerce has given rise to rapid growth in the transport demand for urban-rural logistics. However, the scattered logistics distribution centers and the long transport distance in rural areas maintain high end-to-end costs, severely restricting urban-rural socio-economic development. In response to this double whammy, governments and scholars are committed to promoting the integration of passenger and freight services. Buses are widely accepted as one of the most effective vehicles for sharing vehicle capacity between passengers and parcels [4]. The European Union has taken the lead in proposing the passenger-freight shared mobility (PFSM), that is, realizing the collaborative distribution of small parcels through the existing public transit network to solve the sustainable development challenges of urban-rural public transport systems and the last-kilometer problem of rural logistics. As an active practitioner of this mode, China has launched PFSM in over 1,100 county-level administrative regions, established nearly 50,000 PFSM stops, and opened more than 11,000 PFSM routes as of 2024. It plans to increase the number of PFSM routes to over 20,000 by 2027 [6]. Practice has proven that PFSM not only improves the operational efficiency of urban-rural public transport systems, but also effectively reduces urban-rural logistics costs and carbon emissions, achieving a win-win situation for public transport enterprises and logistics companies [7].

PFSM not only enhances the overall efficiency of urban-rural public transport systems but also significantly alters the operational characteristics, as manifested by the following: (i) PFSM causes cascading fluctuations in bus travel time. Loading and unloading parcels not only increase the dwell time at stops, but also have cascading effects on bus travel time through queuing. These disturbances further affect the punctuality and reliability of passenger services, forming a negative feedback loop between passenger loss and freight demand contraction. (ii) PFSM increases the complexity of bus scheduling and vehicle allocation. Freight loading and unloading requests directly affect the dwell time and bus punctuality at stops to limit the bus schedule. Besides, vehicle type selection and passenger-freight capacity allocation limit determine the maximum transport capacity and constrain vehicle allocation decisions. A bus schedule is required to adapt to differences in passenger/freight capacity between vehicle types, while vehicle allocation is guided by dynamic scheduling requirements. This interaction is further complicated by the temporal and spatial mismatch between passenger and freight demands, resulting in the inapplicability of optimization techniques for traditional urban-rural bus operation and restricting bus capacity utilization under PFSM. (iii) The passenger/freight capacity allocation is a key factor affecting the overall efficiency of the system. An improper design may lead to efficiency losses in both passenger and cargo services. For example, insufficient passenger capacity lowers the turnover capacity of passenger transport and increases the risk of passenger refusal, potentially failing to meet passenger demand. Limited freight capacity restricts loading efficiency, which may result in a waste of freight resources and affect freight profitability.

Aiming at the above issues, this paper focuses on the urban-rural operation under PFSM. Considering the influence of the spatio-temporal distribution characteristics of passenger-freight demand and dynamic stochastic travel time, a bilevel optimization model for multi-type vehicle scheduling and dynamic vehicle capacity allocation is established. The improved jellyfish search algorithm is designed to solve the model. The main contributions of this paper are as follows:

(i) A collaborative optimization methodology combining multi-type vehicle scheduling and dynamic vehicle capacity allocation is introduced for PFSM in urban and rural areas. With the consideration of the dynamic stochastic characteristics of bus travel time and freight loading/unloading, this method achieves the flexible expansion of passenger/freight capacity and differentiated vehicle combinations (passenger-dominated and freight-dominated) to adapt to the spatio-temporal mismatch of urban-rural passenger and freight demands. This contributes to the extension of PFSM-related optimization theories.

(ii) A novel improved jellyfish search algorithm, which is based on the jellyfish search algorithm and introduce tent chaotic map initialization, Levy flight, and differential evolution, is proposed to efficiently to tackle the proposed collaborative optimization problem.

(iii) Impacts of passenger-freight demand, the maximum passenger travel time, and the minimum percentage of passenger capacity to vehicle capacity on optimization schemes are analyzed. The economic and environmental benefits of the proposed method under different scenarios, as well as policy implications and potential applications are presented.

The remainder of this paper is organized as follows. Section 2 presents a literature review of related work. In Section 3, a bilevel optimization model is formulated under this study. Section 4 describes the proposed solution method. Section 5 presents the results and sensitivity analysis of case study, followed by the discussion about policy implications and applications. Finally, Section 6 concludes the paper.

## 2. Literature review

PFSM is essentially a deep integration of passenger transport system and logistics system. This integration not only relies on the guidance of public policies and regulations but also requires close collaboration between traffic participants and logistics suppliers [8]. With the vigorous promotion by government ministries, this practice of utilizing the excess bus capacity to transport short-distance parcels has attracted increasing attention from scholars, giving rise to a series of related studies [9][10]. Early research focused on the utilization in urban public transport networks, while the recent focus has gradually shifted to the dynamic collaborative optimization of spatio-temporal heterogeneous demand in rural scenarios. However, the modelling methods for PFSM vary significantly depending on the application scenario. Urban and rural scenarios differ fundamentally in terms of demand characteristics, core contradictions, and system objectives, leading to modelling differences in the focus, technical approaches, and solution strategies. Therefore, we summarize the existing studies of PFSM in urban and rural areas, respectively.

## 2.1 Urban bus operation under PFSM

Urban public transport system has become the most widely used mean of transport for early PFSM research due to its high frequency and high density [11]. The urban PFSM demand is spatially concentrated (e.g., central business districts and transport hubs) and time-sensitive especially in morning and evening peak hours. The core contradiction lies in the conflict between limited spare vehicle capacity and increasing freight demand. Therefore, most scholars tried to balance on-time delivery with passenger waiting time limit in their modelling to respond to high-density demand through refined resource allocation. For example, Ghilas et al. [12] introduced shuttle buses for transporting passengers and parcels between bus terminals as well as formulated a PFSM bus scheduling model solved by CPLEX solver. However, one of their assumption that parcels and passengers may be forced to wait for a long time at terminals until the departure time, is inconsistent with the facts. Trentini et al. [13] considered the buses transporting parcels from distribution centers to bus stops and proposed an adaptive large neighborhood search, but their study focused on route planning for the last-mile delivery of urban freight fleets. Cheng et al. [14] used mixed-integer linear programming to model the parcel delivery problem for urban public transport and showed that over 60% of parcels are successfully delivered to those destinations within 500 meters of bus stops based on their proposed method. Peng et al. [15] established a two-stage model for PFSM in urban public transport, where the first stage examined the fair matching of passengers and the second stage focused on the parcel delivery by bus. He and Yang [16] proposed a delivery mode for small parcels based on collaboration between passenger transport companies and express delivery companies, with the objective of minimizing total costs by determining express delivery batches, departure times, and distribution routes. Bruno et al. [17] conducted an intensive study on the combinatorial optimization problem of integrating freight services into urban public transport networks by fully considering important constraints (e.g., requirements for freight demand and delivery time windows) and formulating an urban transit network design model under stochastic freight demand. Jonas et al. [18] developed a combinatorial optimization model of vehicle configuration and routing for the modular buses under PFSM, comprehensively considering pick-up and drop-off times, travel time costs, travel distance costs, and fleet size.

Although researchers have achieved valuable results in the studies of urban PFSM, PFSM has certain limitations when applied to urban public transport systems [19]. Compared with rural residents, urban residents have higher population density and travel demand, leading to relatively scarce spare capacity in public transport systems, which may be insufficient to support additional freight services especially in developed cities [20]. If PFSM is mechanically applied to public transport systems in developed cities, it not only fails to effectively offset operating costs but may also has a negative impact on passenger services [21]. This finding has prompted scholars to concentrate on urban-rural bus routes with more idle capacity.

## 2.2 Urban-rural bus operation under PFSM

PFSM in rural areas has been proven to have a greater impact on travel behavior and bus operations compared with that in urban areas [22-24]. Driven by this, scholars are gradually beginning to conduct research and applications around PFSM in the urban-rural context [25].

The urban-rural PFSM demand tends to be distributed asymmetrically in terms of time and space (e.g., seasonal transport of agricultural products and commuting tidal pattern). The key issue of urban-rural PFSM lies in the contradiction between long-distance transport costs and low-frequency demand coverage. Therefore, the existing literature mainly refers to traditional public transport optimization theory or customized public transport operation modes in their modelling strategies, adopting mixed-integer programming to achieve global optimal of resource integration [26]. Wang et al. [27] proposed a mixed-integer programming model to describe urban-rural PFSM in different networks by pre-designing truck routes and some related sets. He et al. [28] addressed the issue of dynamic rural bus scheduling under PFSM. They proposed a hybrid heuristic solve the mixed-integer linear programming model for bus timetable scheduling. Zeng et al. [29] considered the pick-up time window and loading/unloading times for freight services to establish an urban-rural customized bus timetable scheduling model with travel cost minimization. The dynamic programming algorithm and subgradient method were designed to solve their problem. Lin et al. [30] formulated a mixed integer programming for the PFSM multimodal transport system and developed a two-stage solution algorithm to optimize vehicle scheduling and freight allocation.

Although existing research has considerable theoretical value and practical significance in promoting the development of urban-rural PFSM, few studies have considered the cascading fluctuations of freight services on bus travel time, passenger transport punctuality, and service reliability. Moreover, most studies focused on bus route design and timetable scheduling under PFSM. It is also necessary to further develop theories on combinatorial optimization of vehicle scheduling and configuration (allocation). Vehicle scheduling is required to consider the capacity and performance of different vehicle types, while vehicle allocation is optimized according to scheduling requirements. These two problems interact with each other, as a result, only by cooperatively optimizing both vehicle scheduling and allocation can the maximum utilization of public transport capacity be achieved under PFSM. Finally, previous studies generally fixed passenger-freight capacity on each bus, which may lead to waste of vehicle capacity in actual operations and significantly weaken the potential PFSM benefits. Therefore, it is essential to incorporate the passenger-freight capacity allocation into the collaborate optimization of vehicle allocation and scheduling to further improve the overall system efficiency under PFSM.

## 3. Problem formulation

### 3.1 Problem statement

We consider the urban-rural bus operation under PFSM in Fig. 1. The urban-rural bus network is composed of multiple two-way routes with different bus stops. Each urban-rural bus stop can represent a regular bus stop, a distribution center (DC), or an integrated transportation service center (ITSC). Regular bus stops only allow passengers to board or alight, and DCs allow for cargo loading/unloading only while an ITSC can provide passenger and freight services simultaneously. For passenger service, passengers submit a trip request prior to the start of the bus service including the drop-off and pick-up locations, the number of passengers, and the expected departure or arrival time window. For freight service, the parcel information includes the loading and unloading stops, parcel size, and parcel weight.

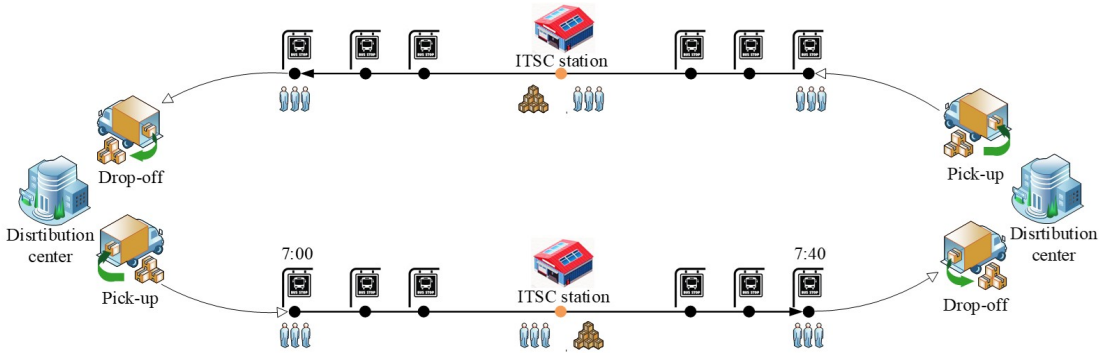

Fig. 1. Illustrative example of a two-way urban-rural bus route under PFSM

Assume that parcels are collected and stored at DCs before being transferred. Before executing the timetable, each bus may travel to the nearest DC to load the parcels, then arrives at its assigned origin station to load passengers and depart on time at a pre-determined time (e.g., 7:00 in Fig. 1). Each bus provides passenger and freight services at the corresponding urban-rural bus stop along the assigned and pre-determined bus routes. Passenger service ends (e.g., 7:40 in Fig. 1) when a bus reaches the final regular bus stop. A bus may visit the nearest DC to unload and load the parcels for its next run, and then arrive at the assigned origin station to serve the opposite direction of the assigned route.

Under the above background, this paper focuses on the collaborative optimization of multi-type bus scheduling and dynamic passenger/freight capacity allocation for urban-rural bus routes under PFSM. Our problem can be formulated as a generalization of the combined optimization problem of vehicle scheduling and allocation. More specifically, the vehicle scheduling problem is required to respond to the dynamic demand at the microscopic level. In order to reduce the impact of freight transport on passenger transport efficiency, the urban-rural bus scheduling problem under PFSM generally chooses to minimize the daily passenger travel time to determine the daily schedule for each bus under a set of constraint limitations (e.g., the bus fleet size, potential bus types, potential passenger and freight capacity configurations). For the vehicle allocation problem, minimizing the total cost or the total profit is the major objective function to determine the optimal fleet size of each vehicle type and the optimal passenger and freight capacities on each bus at the macroscopic level. In other words, the vehicle scheduling problem can receive its input from the vehicle allocation problem, which could include parameters, constraints, or even decision variables of the vehicle allocation problem, while the vehicle allocation problem should consider the anticipated behavior of the outputs from the vehicle scheduling problem as well. Therefore, our proposed problem can be regarded as a bilevel programming problem [31]. We propose a bilevel optimization model to formulate multi-type bus scheduling and dynamic passenger/freight capacity allocation for urban-rural bus routes. In the upper level model, the vehicle type and the passenger/freight capacity allocation on each bus are optimized to constitute some feasible boundaries of the lower level problem. The lower level problem is the bus scheduling problem to assign the bus runs to a set of available vehicles and then determine the total passenger travel time during the daily operation, influencing the passenger travel cost and the revenue in the upper level. As a result, the solution to the lower level problem affects the upper level objective function value and the bus allocation.

### 3.2 Assumptions and notations

The basic assumptions of the study are made as follows: (i) origin-destination (OD) passenger-freight trip information is available before the bus operation; (ii) each urban-rural bus has limited passenger/freight capacity based on the vehicle type; (iii) each urban-rural bus starts and ends its trip at the same depot, and provides passenger and freight services simultaneously in ITSCs; (iv) we only consider small light-weight parcels for freight service in this paper; (v) the bus departure timetable of the bus routes is known.

Let $R$, $N$, and $K$ be the sets of bus runs, vehicle types, and bus numbers, respectively, where $r \in R$, $n \in N$, $k \in K$. $S$ is the set of bus stops, where $\forall i \in S$. We introduce binary variables $\delta_i$, $x_k^r$, and $y_k^n$. $\delta_i$ equals 1 if stop $i$ allows for cargo loading/unloading; and 0 otherwise. $x_k^r$ equals 1 if bus $k$ serves run $r$; and 0 otherwise. $y_k^n$ equals 1 if the type of bus $k$ is type $n$; and 0 otherwise. $\lambda_k$ is the percentage of passenger capacity to vehicle capacity on bus $k$. Other notations of the problem are summarized in Appendix A1.

### 3.3 The upper level model

The upper level problem is the bus allocation problem. Under PFSM, we consider multiple vehicle types and dynamic passenger-freight capacity allocation in this paper, which can directly affect the system capacity supply structure in the following respects: (1) the difference in the purchasing cost of multiple vehicle types forms the capital expenditure constraint; (2) the passenger-freight capacity allocation determines the transport efficiency; and (3) the constraint of marginal cost is based on the correlation between the road tolls and the vehicle types. The objective is to maximize the operating profit and obtain the Pareto optimization by constructing the purchasing cost $C_{\text{fix}}$, the running cost $C_{\text{km}}$, the dwelling cost $C_{\text{dwell}}$, and the toll cost $C_{\text{toll}}$, and estimating the passenger-freight revenue. Under profit maximization, the decision variables are vehicle type $y_k^n$ and the percentage of passenger capacity to vehicle capacity $\lambda_k$ for bus $k$, which also define the hard constraints of vehicle scheduling in the lower level.

#### 3.3.1 Toll cost

Whether urban-rural buses are required to pay for road tolls depends on the use of toll roads and the local transport policy. If a bus route does not pass through a toll road or toll-free policy is implemented, $C_{\text{toll}} = 0$. In some countries (e.g., China), the road toll setting is required based on the local policy, passenger capacity (vehicle seat number), and mileage rate. $C_{\text{toll}}$ can then be described as:

$$C_{\text{toll}} \in \left\{ 0, \sum_{r \in R} \sum_{n \in N} \sum_{k \in K} x_k^r y_k^n \Delta_{\text{toll}} \left( L_{\text{toll}}, \beta_k^n \right) \right\} \tag{1}$$

where $L_{\text{toll}}$ is the toll mileage of a bus route in km; $\Delta_{\text{toll}} \left( L_{\text{toll}}, \beta_k^n \right)$ is the road toll function of a bus; $\beta_k^n$ is the seat number of bus $k$ with type $n$.

For vehicles with more than seven seats, the tollway calculations are divided into four classes based on their passenger capacity on most Chinese expressways. $\Delta_{\text{toll}} \left( L_{\text{toll}}, \beta_k^n \right)$ can be calculated as:

$$\Delta_{\text{toll}} \left( L_{\text{toll}}, \beta_k^n \right) = \begin{cases} \alpha_1 L_{\text{toll1}}, & \beta_k^n \leq \beta_1 \\ \alpha_2 L_{\text{toll}}, & \beta_1 < \beta_k^n \leq \beta_2 \\ \alpha_3 L_{\text{toll}}, & \beta_2 < \beta_k^n \leq \beta_3 \\ \alpha_4 L_{\text{toll}}, & \beta_k^n > \beta_3 \end{cases} \tag{2}$$

$$\beta_k^n = \left\lfloor \frac{1}{v_{\text{p}}} \lambda_k V_n \right\rfloor \tag{3}$$

where $\bar{v}_p$ denotes the average spatial volume of a seat in a bus in m³; $\alpha_i$ is the toll rate for class $i$ ($i$ = 1, 2, 3, 4); $\beta_j$ is the seat number threshold for different classes ($j$ = 1, 2, 3); $V_n$ is the vehicle capacity of bus type $n$ in m³; $\lfloor \cdot \rfloor$ is the floor function.

### 3.3.2 Dwelling cost

The total dwelling cost of bus stops depends on the loading and uploading amounts of passengers/parcels at each stop. It can be calculated as:

$$C_{\text{dwell}} = \frac{\phi_{\text{dwell}}}{3600} \sum_{r \in R} \sum_{m \in S} \left( \Delta^{1}_{r,m} + \Delta^{-1}_{r,m} \right) \tag{4}$$

where $\phi_{\text{dwell}}$ is the value of dwell time; $\Delta^{1}_{r,m}$ and $\Delta^{-1}_{r,m}$ are the dwell times of a bus serving the $r$th run in the forward direction and the opposite direction at stop $m$ in seconds, respectively.

In order to better describe the route directions, let $|i|$ denotes the stop sequence number of stop $i$ in the bus network, then the route direction can be expressed by:

$$u(i,j) = \frac{|j| - |i|}{\||j| - |i|\|} \tag{5}$$

where $u(i,j) = 1$ if stop $i$ locates upstream of stop $j$ in the stop sequence; and -1 otherwise.

The dwell times of a bus serving the $r$th run at stop $m$ is the maximum value among the loading and uploading times of passengers/parcels, which can be described as:

$$\Delta^{u(i,j)}_{r,m} = \max \left\{ \Delta_p \sum_{i \in S} d^{r}_{i,m}, \Delta_p \min \left\{ \sum_{j \in S} d^{r}_{m,j}, P^{re}_{r,m} \right\}, \Delta_f \left( q^{r}_{i,m} + q^{r}_{m,j} \right) \delta_i \delta_j \right\}, i < m < j \tag{6}$$

$$P^{re}_{r,i} = \sum_{n \in N} \sum_{k \in K} \beta_{k,n} X_{k,r} Y_{k,n} - \sum_{j=1}^{i-1} \left( o^{u}_{r,j} - d^{u}_{r,j} \right) + d^{u}_{r,i}, \forall r \in R, i \in S, i \geq 2 \tag{7}$$

where $\Delta_p$ is the average delay of a boarding/alighting passenger in seconds; $\Delta_f$ is the average delay of a loading/unloading parcel in seconds; $d^{r}_{i,j}$ and $q^{r}_{i,j}$ are the number of passengers and the amount of cargo transported from stop $i$ to stop $j$ in the $r$th run, respectively; $o^{u}_{ri}$ and $d^{u}_{ri}$ are the number of passengers boarding and alighting at stop $i$ on the bus serving the $r$th run, respectively; $\delta_i$ is a binary variable that is 1 if stop $i$ is allowed for cargo loading/unloading; and 0 otherwise. $P^{re}_{r,i}$ is the number of seats remaining when the bus serving the $r$th run arrives at stop $i$, which is related to the vehicle type, the percentage of passenger capacity to vehicle capacity, the number of passengers alighting at stop *i*, and the numbers of passengers boarding and alighting at the upstream stops.

### 3.3.3 Purchasing and running costs

The running cost is mainly related to the running distance and the vehicle type for those corresponding bus runs. Therefore, the total running cost $C_{km}$ is formulated as follows:

$$C_{km} = \sum_{r \in R} \sum_{n \in N} \sum_{k \in K} L_r \phi_n^{km} x_k^r y_k^n \tag{8}$$

where $L_r$ is the running distance of the *r*th run in km; $\phi_n^{km}$ is the running cost per kilometer of bus type *n*.

The purchasing cost of the bus fleet is dependent on the bus type and the corresponding number of vehicles, i.e.:

$$C_{fix} = \sum_{n \in N} \sum_{k \in K} \phi_n^{buy} y_k^n \tag{9}$$

where $\phi_n^{buy}$ is the purchasing cost of bus type *n*, including vehicle depreciation, repair, and maintenance, and driver wages.

### 3.3.4 Revenue

Transit passenger fare revenue $E_u$ correlates to ridership, passenger miles travelled, and vehicle type, which can be calculated as follows:

$$E_u = \begin{cases} \sum_{r \in R} \sum_{i \in S} \sum_{j \in S} d_{i,j}^r \left( l_{i,j} \sum_{n \in N} \sum_{k \in K} \eta_n^u x_k^r y_k^n + \gamma_u \right), & l_{i,j} \leq L_u \\ \sum_{r \in R} \sum_{i \in S} \sum_{j \in S} d_{i,j}^r \left( (l_{i,j} - L_u) \sum_{n \in N} \sum_{k \in K} \eta_n^u x_k^r y_k^n + \gamma_u \right), & l_{i,j} > L_u \end{cases} \tag{10}$$

where $\gamma_u$ is the base fare of passenger service for the first $L_u$ kilometers and $\eta_n^u$ is charged thereafter per kilometer of bus type *n*.

Currently, there are few specific studies showing that cargo shipping costs are related to the vehicle types. We assume that the revenue of freight service $E_f$ is only related to the volume of freight and traveled distance, $E_f$ can then be calculated as:

$$E_{\mathrm{f}} = \begin{cases} \sum_{r \in R} \sum_{i \in S} \sum_{j \in S} q_{i,j}^r \left( \eta_{\mathrm{f}} l_{i,j} + \gamma_{\mathrm{f}} \right), l_{i,j} \leq L_{\mathrm{f}} \\ \sum_{r \in R} \sum_{i \in S} \sum_{j \in S} q_{i,j}^r \left( \eta_{\mathrm{f}} \left( l_{i,j} - L_{\mathrm{f}} \right) + \gamma_{\mathrm{f}} \right), l_{i,j} > L_{\mathrm{f}} \end{cases} \quad (11)$$

where $\gamma_{\mathrm{f}}$ is the base fare of freight service for the first $L_{\mathrm{f}}$ kilometers and $\eta_{\mathrm{f}}$ is charged thereafter per kilometer.

### 3.3.5 Model formulation

Based on the above definitions, the mathematical model of the upper level problem is shown as follows.

$$\max_{\lambda, x, y} Z = E_{\mathrm{u}} + E_{\mathrm{f}} - C_{\mathrm{km}} - C_{\mathrm{idle}} - C_{\mathrm{fix}} - C_{\mathrm{toll}} \quad (12)$$

s.t.

$$\sum_{i \in S} \left( o_{r,i}^{\mathrm{u}} - d_{r,i}^{\mathrm{u}} \right) \leq \sum_{n \in N} \sum_{k \in K} \lambda_k \beta_k^n x_k^r y_k^n, \forall r \in R, i \geq 1 \quad (13)$$

$$\sum_{i \in S} \left( o_{r,i}^{\mathrm{f}} - d_{r,i}^{\mathrm{f}} \right) \leq \sum_{n \in N} \sum_{k \in K} (1 - \lambda_k) V_n x_k^r y_k^n, \forall r \in R, i \geq 1 \quad (14)$$

$$\sum_{n \in N} y_k^n = 1, \forall k \in K \quad (15)$$

$$\lambda_{\min} \leq \lambda_k \leq 1, \forall k \in K \quad (16)$$

$$\delta_i \in \{0,1\}, \forall i \in S \quad (17)$$

$$x_k^r \in \{0,1\}, \forall r \in R, k \in K \quad (18)$$

$$y_k^n \in \{0,1\}, \forall n \in N, k \in K \quad (19)$$

Objective (12) is to maximize the operating profit $Z$.

Constraints (13) and (14) guarantee that the number of passengers and the amount of cargo transported by the bus serving the $r$th run is not greater than the passenger and freight capacities, respectively. Constraint (15) ensures that the bus type of the same bus is unique. Constraint (16) sets a lower bound for the percentage of passenger capacity to vehicle capacity due to the minimum passenger capacity requirement in practice. Constraints (17)-(19) define the domains of the decision variables.

### 3.4 The lower level model

The lower level problem is the bus scheduling problem. We focus on the collaborative optimization of passenger and cargo transports spatially and temporally at the operation level. We consider the following problems: (1) cargo loading/unloading triggers fluctuations in stop dwell time; (2) differences in the spatial and temporal distributions of passenger and freight demands lead to the conflicts in vehicle scheduling; and (3) the variation in travel time affects the bus reliability. Since passenger transport has higher requirements for travel timeliness, our objective is to minimize the total passenger travel time including total in-vehicle travel time, total passenger detention time, total dwell time, and total passenger waiting time in the lower level. Besides, the stochasticity of the travel time is quantified by constructing the improved Bureau of Public Road (BPR) function. As the decision variables in the lower level, the optimal bus scheduling can influence the passenger travel cost and the revenue in the upper level. As a result, the solution to the lower level problem affects the upper level objective function value and the bus allocation to achieve better bus service efficiency.

### 3.4.1 Total in-vehicle travel time

In urban-rural bus operation, in-vehicle travel time cost is a key element that affects system performance. The calculation of this time cost term is required to comprehensively consider the spatio-temporal characteristics of the vehicle operating state and the dynamic characteristics of the passenger demand distribution. Therefore, we calculate in-vehicle travel time incorporating dynamic stochastic travel time reliability. The total in-vehicle travel time $T_{\text{cruise}}$ can be expressed as:

$$T_{\text{cruise}} = \sum_{r \in R} \sum_{i \in S} \sum_{j \in S} \sum_{k \in K} d_{i,j}^r E\left[T_{i,j}^r\right] x_k^r, i \neq j \tag{20}$$

where $E\left[T_{i,j}^r\right]$ is the expected in-vehicle travel time between $(i, j)$ of the bus serving the $r$th run.

In order to capture the time-varying characteristics of urban-rural road networks, this paper adopts the improved BPR function [32] to estimate the in-vehicle travel time between stops:

$$T_{i,j}^r = a_{i,j}^r \left[1 + \beta \left(\frac{Q_{i,j}^r}{C_{i,j}^r}\right)^z\right] + \varepsilon_{i,j}^r \tag{21}$$

where $T_{i,j}^r$, $a_{i,j}^r$, $Q_{i,j}^r$, and $C_{i,j}^r$, are in-vehicle travel time, free-flow travel time, demand volume, and capacity of the bus serving the $r$th run between $(i, j)$, respectively; $\varepsilon_{i,j}^r$ is the stochastic perturbation term under the normal distribution $N(0, \sigma^2)$; $\beta$ and $z$ are the calibration parameters. The estimation takes into account both the deterministic congestion effect and the stochastic fluctuation factor, providing a more realistic representation of real-world traffic conditions compared with the traditional BPR function.

Aiming at the travel time uncertainty, this paper introduces the time budget theory to formulate the reliability constraint, which enables the bus scheduling scheme to cope with both normal traffic congestion and unexpected delays, providing a strict mathematical basis for subsequent multi-objective optimization. We define $\hat{T}_{i,j}^{r}$ as the time budget value that satisfies the confidence level, which is the sum of the expected in-vehicle travel time between stops and the buffer time for departure or arrival delays. It can be described as:

$$\hat{T}_{i,j}^{r} = \inf\left\{\tau \in \mathbb{R}^{+} \mid P\left(T_{i,j}^{r} \leq \tau\right) \geq \gamma\right\} \tag{22}$$

When $T_{i,j}^{r}$ follows a normal distribution, its analytical expression is obtained as follows:

$$\hat{T}_{i,j}^{r} = \mu_{i,j}^{r} + \sigma_{i,j}^{r} \cdot z_{\gamma} \tag{23}$$

where $\tau$ is the potential travel time candidates; $\mu_{i,j}^{r} = E\left[T_{i,j}^{r}\right]$ and $\sigma_{i,j}^{r} = \sqrt{Var\left(T_{i,j}^{r}\right)}$; $\Phi_{\gamma}$ is the $\gamma$-quantile of a standard normal distribution, and its correction $\Phi_{\gamma}'$ can be obtained by the Cornish-Fisher expansion considering higher order moments of returns [33], as shown in Equation (24):

$$\Phi_{\gamma}' = \Phi_{\gamma} + \frac{\tilde{S}}{6}\left(\Phi_{\gamma}^{2} - 1\right) + \frac{\tilde{K}}{24}\left(\Phi_{\gamma}^{3} - 3\Phi_{\gamma}\right) - \frac{\tilde{S}^{2}}{36}\left(2\Phi_{\gamma}^{3} - 5\Phi_{\gamma}\right) \tag{24}$$

where $\tilde{S}$ and $\tilde{K}$ are the skewness and kurtosis of the in-vehicle travel time, respectively. By substituting $\Phi_{\gamma}'$ into Equation (23), the accuracy of time budget estimation can be improved under the condition of non-normal distributions.

To ensure service reliability, we introduce an indicator for in-vehicle travel time reliability. $R_{i,j}^{r}$, the probability that the actual in-vehicle travel time from stops $i$ to $j$ of the bus serving the $r$th run is not greater than the travel time budget, is applied:

$$R_{i,j}^{r} = P\left\{T_{i,j}^{r} \leq \hat{T}_{i,j}^{r}\right\} \geq \gamma \tag{25}$$

Based on the law of large numbers and the central limit theorem, the 85% confidence level corresponds to the Z-value of 1.036 in a standard normal distribution, which can effectively cover the main probability intervals of random perturbations (68.27% in the range of $\pm 1\sigma$, and 86.64% the range of $\pm 1.5\sigma$). Corrected by Cornish-Fisher expansion, this confidence level is compatible with the right-skewed distribution characteristics, commonly found in the actual traffic flow (e.g., positively skewed distribution of travel time). Therefore, this paper sets $\gamma = 0.85$ as the reliability threshold, i.e., 85% of the actual in-vehicle travel time is guaranteed to satisfy the time budget constraint. The value is selected to balance the service reliability and operational efficiency, which not only avoids the over-conservative resource allocation, but also avoids the risk of demand loss caused by the lack of reliability, and meets the dual needs for punctuality and profitability in the transit system.

### 3.4.2 Total passenger detention time

When the peak passenger flow is high and bus capacity is insufficient, some passengers may not be able to board the bus. Assume that there are no secondary passenger detentions, the total passenger detention time $T_{\text{det}}$ can be calculated as:

$$T_{\text{det}} = \sum_{r \in R} \sum_{i \in S} \left( t_{r+1,i} - t'_{r,i} \right) P_{r,i}^{\text{det}} \tag{26}$$

$$P_{r,i}^{\text{det}} = \max \left\{ 0, P_{r,i}^{\text{wait}} - P_{r,i}^{\text{re}} \right\}, \forall r \in R, i \in S \tag{27}$$

where $P_{r,i}^{\text{det}}$ is the number of stranded passengers after the $r$th run leaves stop $i$, depending on the difference between the number of passengers waiting for the $r$th run ($P_{r,i}^{\text{wait}}$) and the number of seats remaining on the bus ($P_{r,i}^{\text{re}}$).

The number of passengers waiting for the $r$th run at stop $i$ is the sum of the cumulative number of arrivals during the departure interval and the number of passengers stranded on the previous run, i.e.,

$$P_{r,i}^{\text{wait}} = \int_{t'_{r-1,i}}^{t_{r,i}} D_i(t) dt + P_{r-1,i}^{\text{det}}, \forall r \in R, i \in S, r \geq 2 \tag{28}$$

The number of seats remaining on the bus serving the $r$th run and arriving at stop $i$ is related to the bus type, the percentage of passenger capacity to vehicle capacity, the number of passengers alighting at stop $i$, and the number of passengers boarding and alighting at upstream stops, i.e.,

$$P_{r,i}^{\text{re}} = \sum_{n \in N} \sum_{k \in K} \beta_k^n x_k^r y_k^n - \sum_{j=1}^{i-1} \left( o_{r,j}^{\text{u}} - d_{r,j}^{\text{u}} \right) + d_{r,i}^{\text{u}}, \forall r \in R, i \in S \text{ 且 } i \geq 2 \tag{29}$$

### 3.4.3 Total dwell time and total passenger waiting time

The total dwell time $T_{\text{dwell}}$ is mainly affected by the loading and uploading passengers/parcels at intermediate stops, which is calculated as follows:

$$T_{\text{dwell}} = \frac{1}{3600} \sum_{r \in R} \sum_{(i,j) \in S} \sum_{k \in K} d_{i,j}^r x_k^r \left( \sum_{\substack{i<n<j \\ u(i,j)=u(m,n)}} \Delta_{r,n}^{u(i,j)} + \sum_{\substack{i<m<j \\ u(i,j)=u(m,n)}} \Delta_{r,m}^{u(i,j)} \right) \tag{30}$$

where $\sum_{\substack{n \in (i,j) \\ u(i,j)=u(m,n)}} \Delta_{r,n}^{u(i,j)}$ denotes the delay caused by passengers alighting at intermediate stops;

$\sum_{\substack{m \in (i,j) \\ u(i,j)=u(m,n)}} \Delta_{r,m}^{u(i,j)}$ denotes the delay caused by passengers boarding at intermediate stops.

The total passenger waiting time $T_{\text{wait}}$ is dependent on the time interval between the adjacent runs to the same stop and the passenger arrival function, which is calculated as follows:

$$T_{\text{wait}} = \sum_{r \in R} \sum_{i \in S} \int_{t'_{r-1,i}}^{t_{r,i}} D_i(t) dt \tag{31}$$

where $D_i(t)$ is the cumulative passenger arrival function at stop $i$; $t'_{r-1,i}$ is the time when the (r-1)-th run departs from stop $i$; $t_{r,i}$ is the time when the $r$th run arrives at stop $i$.

### 3.4.4 Model formulation

Based on the above definitions, the mathematical model of the lower level problem is shown as follows.

$$\min_x T = T_{\text{cruise}} + T_{\text{dwell}} + T_{\text{wait}} + T_{\text{det}} \tag{32}$$

s.t.

$$o_{r,i}^{\text{u}} = \sum_{j \in S} d_{i,j}^r, \forall r \in R, i \in S, i < j \tag{33}$$

$$d_{r,j}^{\text{u}} = \sum_{i \in S} d_{i,j}^r, \forall r \in R, j \in S, i < j \tag{34}$$

$$o_{r,i}^{\text{f}} = \sum_{j \in S} q_{i,j}^r, \forall r \in R, i \in S, i < j \tag{35}$$

$$d_{r,j}^{\text{f}} = \sum_{i \in S} q_{i,j}^r, \forall r \in R, j \in S, i < j \tag{36}$$

$$\frac{T}{\sum_{r \in R} \sum_{i \in S} \sum_{j \in S} d_{i,j}^r} \leq T_{\max} \tag{37}$$

$$\sum_{r \in R} \sum_{i \in S} \sum_{j \in S} d_{i,j}^r \leq \sum_{r \in R} \sum_{n \in N} \sum_{k \in K} \lambda_k x_k^r y_k^n V_n \tag{38}$$

$$\sum_{r \in R} \sum_{i \in S} \sum_{j \in S} q_{i,j}^r \leq \sum_{r \in R} \sum_{n \in N} \sum_{k \in K} (1 - \lambda_k) x_k^r y_k^n V_n \tag{39}$$

$$\sum_{k \in K} x_k^r = 1, \forall r \in R \tag{40}$$

Objective (32) is to minimize the total passenger travel time $T$.

Constraints (33) and (34) are passenger flow conservation constraints at bus stops. Constraints

(35) and (36) are freight flow conservation constraints at bus stops. Constraint (37) determines the upper bound of the average passenger travel time $T_{max}$. Constraints (38) and (39) guarantee that the vehicle capacity can cover the passenger-freight demand. Constraint (40) ensures that each run can only be served by one bus.

### 3.5 Model transformation for bilevel multi-objective optimization

The proposed model is bilevel multi-objective program with nonlinear constraints. To efficiently solve for solutions to this model, we use the entropy weight method (EWM) to combine the upper-level objective (max $Z$, i.e., profit maximization) and the lower-level objective (min $T$, i.e., total travel time minimization) into a single objective function $F$, transforming the original model into a single-level single-objective program to maximize $F$:

$$\max F = w_1 Z - w_2 T \tag{41}$$

where $F$ can be regarded as the system benefit; $w_1$ and $w_2$ are the weighting factors for $Z$ and $T$, respectively.

EWM can measure the uncertainty in decision-making by calculating the information entropy of each index, which can avoid the anthropogenic interference and make the evaluation more in line with the reality. It has been widely used in transportation optimization problems [34][35]. In order to eliminate the influence of the dimension of the evaluation criteria, the minimum-maximum normalization of the total travel time and total operating profit is applied. The information entropy of the $j$th evaluation index $E_j$ ($j = 1$ for $T$ and $j = 2$ for $Z$) can be calculated as:

$$E_j = -\eta^* \sum_{i=1}^{O} p_{i,j} \ln p_{i,j}, j = 1, 2 \tag{42}$$

where $O$ is the number of alternatives; $\eta^* = 1/\ln n^*$; $n^*$ is the number of evaluation indices, $n^* = 2$ in this paper.

$p_{i,j}$ is the weight of the $i$th alternative based on $j$th evaluation index, which can be calculated as follows:

$$p_{i,j} = \frac{r_{i,j}}{\sum_{i=1}^{O} r_{i,j}} \tag{43}$$

where $r_{i,j}$ is the performance value of the $i$th alternative based on $j$th evaluation index.

Based on the value of the information entropy, the weighting factor of the $j$th index $w_j$ can

be calculated by:

$$w_j = \frac{(1-E_j)}{\sum_{j=1}^{3}(1-E_j)} \tag{44}$$

## 4. Solution method

To efficiently solve for solutions to this model, a heuristic or metaheuristic is required. Most heuristics have some common problems when dealing with complex optimization problems: they may converge to the local optimal solution too early, or the convergence speed is slow. In order to overcome the above problems, we propose an improved jellyfish search algorithm (IJS) to solve our model. The schematic flowchart of IJS is shown in Fig. 2, and the specific steps are described in Sections 4.1 and 4.2.

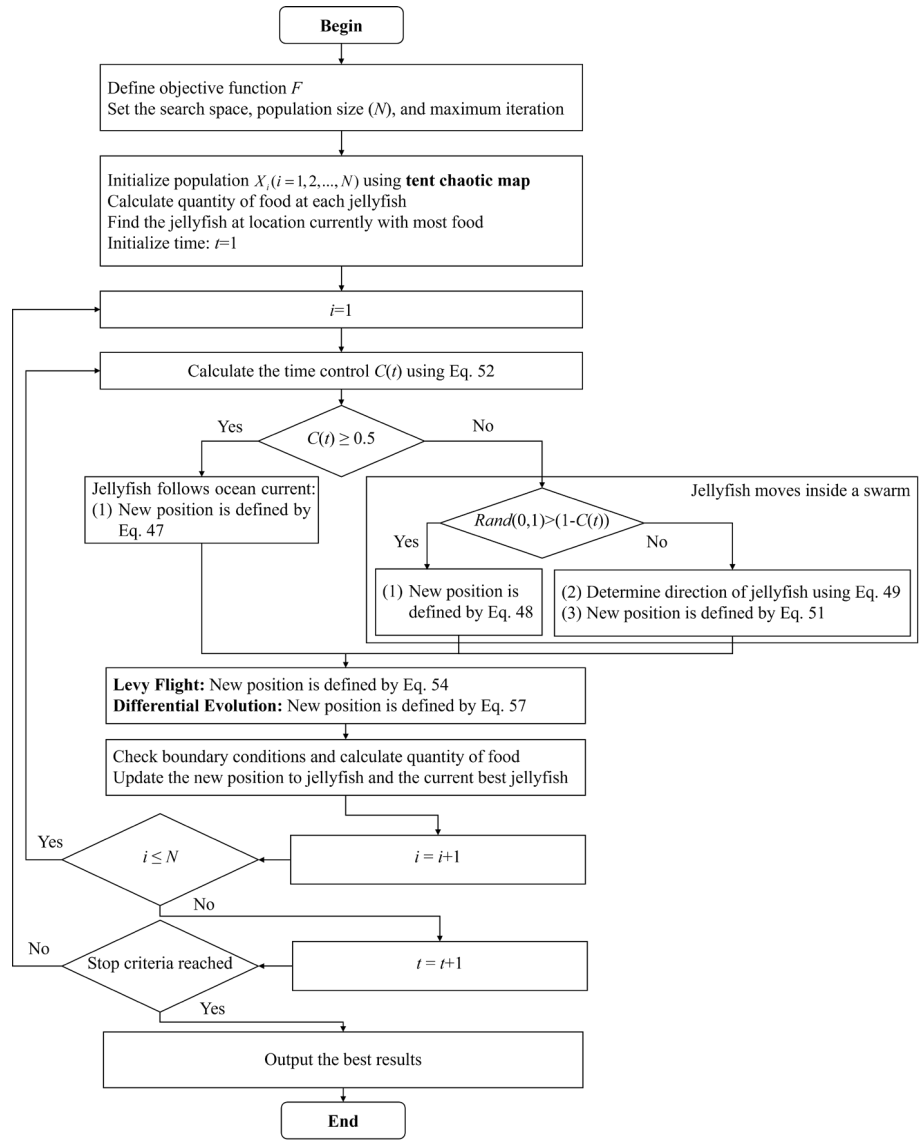

Fig. 2. Schematic flowchart of IJS

## 4.1 Solution representation

In this paper, integer encoding is used to represent the scheme for bus allocation and scheduling. Let $K^*$ and $R^*$ be the total fleet size and the daily number of bus runs, respectively. A solution is a scheme for bus allocation and scheduling $X=\{x, y, \lambda\}$ containing three parts: the first part is the bus scheduling scheme $x$ with dimension of $1\times R^*$; the second is the bus type allocation scheme $y$ with dimension of $1\times K^*$; and the third is the bus passenger capacity allocation scheme $\lambda$ with dimension of $1\times K^*$. The total dimension of solution $X$ is $D=1\times(R^*+2K^*)$. Fig. 3 shows an example of a scheme with $K^*=6$ and $R^*=12$.

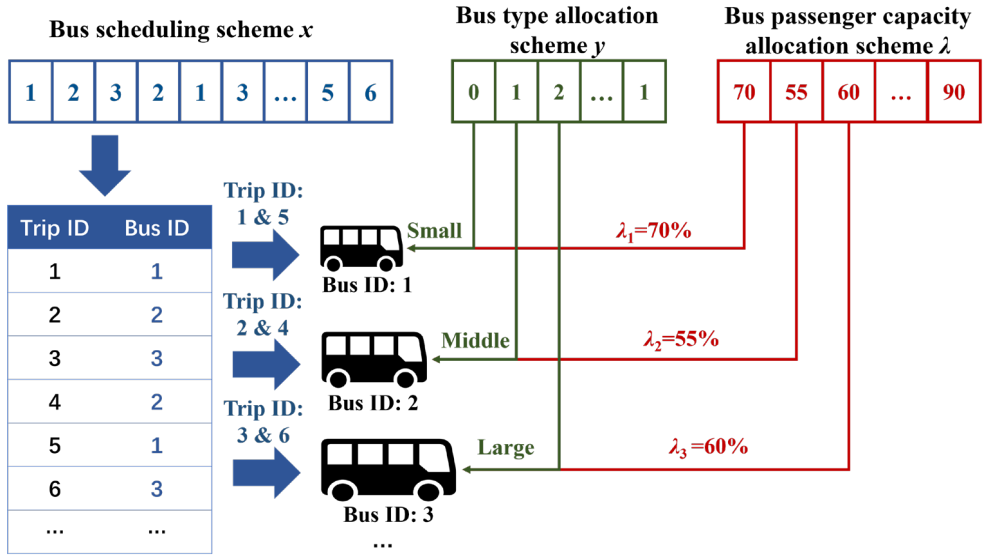

Fig. 3. An example of an encoding scheme $X$ with $K^*=6$ and $R^*=12$

In Fig. 3, $X=\{[1, 2, 3, 2, 1, 3, 4, 5, 6, 4, 5, 6], [0, 1, 2, 2, 0, 1], [70, 55, 60, 80, 100, 90]\}$. The bus scheduling scheme $x$ is [1, 2, 3, 2, 1, 3, 4, 5, 6, 4, 5, 6], where each element denotes the bus number serving the corresponding run. The bus type allocation scheme $y$ is [0, 1, 2, 2, 0, 1], where each element denotes the vehicle type of each bus (0, 1, and 2 represent small, medium, and large, respectively). The bus passenger capacity allocation scheme $\lambda$ is [70, 55, 60, 80, 100, 90], where each element denotes the percentage of passenger capacity to vehicle capacity on each bus. For example, No. 3 bus is a large-size vehicle with 60% vehicle capacity for passenger service to serves the 3rd and 6th runs.

## 4.2 Improved Jellyfish Search

The proposed IJS is an enhanced version of the jellyfish search algorithm (JS) [36] proposed by Cho and Truong (2021) to improve their framework of searching for a global optimal solution. To mimic the food search behavior of jellyfish in the ocean, a scheme for bus allocation and scheduling $X_i$ corresponds to the location of jellyfish *i*. The fitness (i.e., food quantity) of the location is quantified by objective function $F(X_i)$. The IJS is based on three idealized rules: (i) jellyfish either follow the ocean current or move within the swarm. The switch between these two movements is handled by a time control mechanism; (ii) jellyfish are attracted to locations with a higher food quantity; (iii) The food quantity found by jellyfish is determined by the corresponding location.

### 4.2.1 Population initialization

Different from the logistic chaotic map used by the JS to initialize jellyfish population, we adopt the tent chaotic map in the proposed IJS, which has better traversal, greater uniformity and faster iteration speed. Tent chaotic map is more suitable for generating chaotic sequences in our problem [37]. The population initialization can be expressed as follows:

$$X_{i+1} = \begin{cases} \mu_{tent} X_i, & X_i < 0.5 \\ \mu_{tent}(1-X_i), & X_i \geq 0.5 \end{cases} \tag{46}$$

where $X_i$ is the initial location of the *i*th jellyfish, which is generated by the tent chaotic map; $\mu_{tent}$ is a parameter that controls the distribution pattern of the tent chaotic sequence.

### 4.2.2 Ocean current

Jellyfish are attracted to ocean currents because they contain large amounts of nutrients. The direction of the current is determined by calculating the average vector of the vectors from the location of each jellyfish to that of the jellyfish with the best objective function in the swarm. Equation (47) updates the location of each jellyfish within the current:

$$X_i^{t+1} = X_i^t + rand(0,1) \times (X_{best}^t - \beta_d \times rand(0,1) \times \bar{\mu}) \tag{47}$$

where $X_i^t$ is the location of the *i*th jellyfish in the *t*th iteration; $X_{best}^t$ is the best location in the current jellyfish population; $\beta_d$ is the distribution coefficient; $\bar{\mu}$ is the mean location of all jellyfish in the population.

### 4.2.3 Jellyfish swarm

When the jellyfish swarm initially forms, most of the jellyfish engage in group movement, i.e. passive movement. New locations are generated as follows:

$$X_i^{t+1} = X_i^t + \gamma' \times rand(0,1) \times (U_b - L_b) \tag{48}$$

where $\gamma'$ is a motion coefficient related to the length of motion around the locations of the jellyfish; $U_b$ and $L_b$ are the upper bound and lower bound of search spaces, respectively.

Over time, the jellyfish take the initiative to move closer to their companions with more food found. Their directions of movement are shown in Equation (49), i.e. active movement. The generation of new locations is shown in Equations (50) and (51). A schematic diagram of jellyfish movement is shown in Fig. 4.

$$\overrightarrow{Direction} = \begin{cases} X_j^t - X_i^t, F(X_j^t) \geq F(X_i^t) \\ X_i^t - X_j^t, F(X_j^t) < F(X_i^t) \end{cases} \quad (49)$$

$$\overrightarrow{Step} = rand(0,1) \times \overrightarrow{Direction} \quad (50)$$

$$X_i^{t+1} = X_i^t + \overrightarrow{Step} \quad (51)$$

where $X_i^t$ is the location of the $i$th jellyfish in the $t$th iteration; $X_j^t$ is a jellyfish randomly selected from the population.

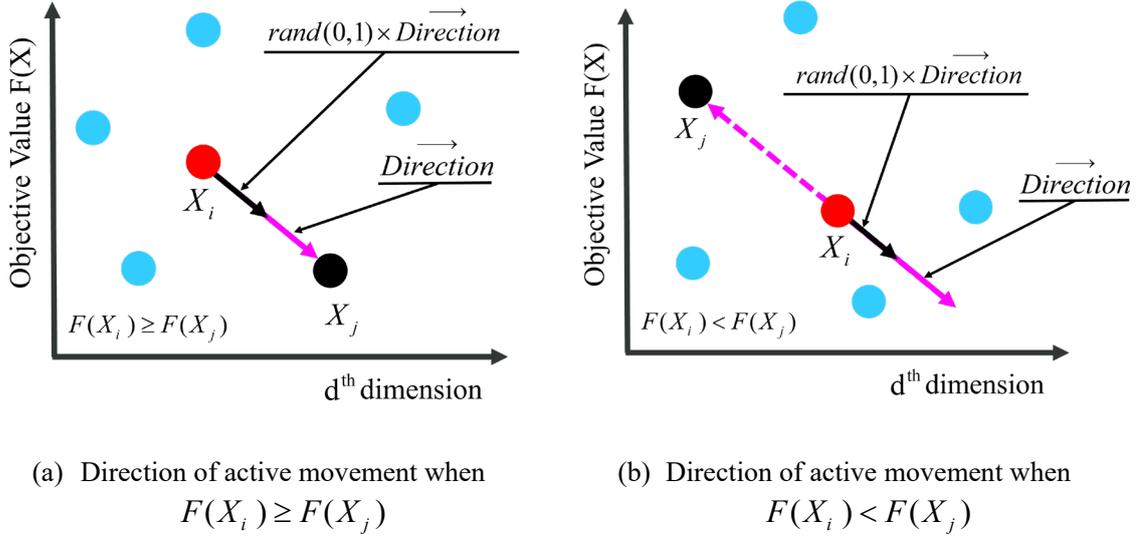

(a) Direction of active movement when $F(X_i) \geq F(X_j)$

(b) Direction of active movement when $F(X_i) < F(X_j)$

Fig. 4. Direction of movement in a swarm

### 4.2.4 Time control mechanism

Over time, jellyfish gather together in ocean currents to form a swarm. When the current is changed, the jellyfish enter another current and create another swarm. Passive movement begins preferentially within a jellyfish swarm. As time passes, individual jellyfish tend to move actively. In order to regulate passive and active movements, the IJS introduces a time control mechanism that regulates the switch. This consists of a time control function $C(t)$ and a threshold constant $C_0$, where the control function is:

$$C(t) = \left|(1-\frac{t}{Max_{it}}) \times (2 \times rand(0,1)-1)\right| \qquad (52)$$

where $C(t)$ is a random value that fluctuates in the interval [0, 1] as the number of iterations increases, $C_0 = 0.5$; $Max_{it}$ is the maximum number of iterations. When $C(t) \geq C_0$, the jellyfish follow the ocean current for movement; conversely, the jellyfish move only within the swarm, i.e., active and passive movement. Fig. 5 simulates the current phase, jellyfish swarm phase, and time control mechanism.

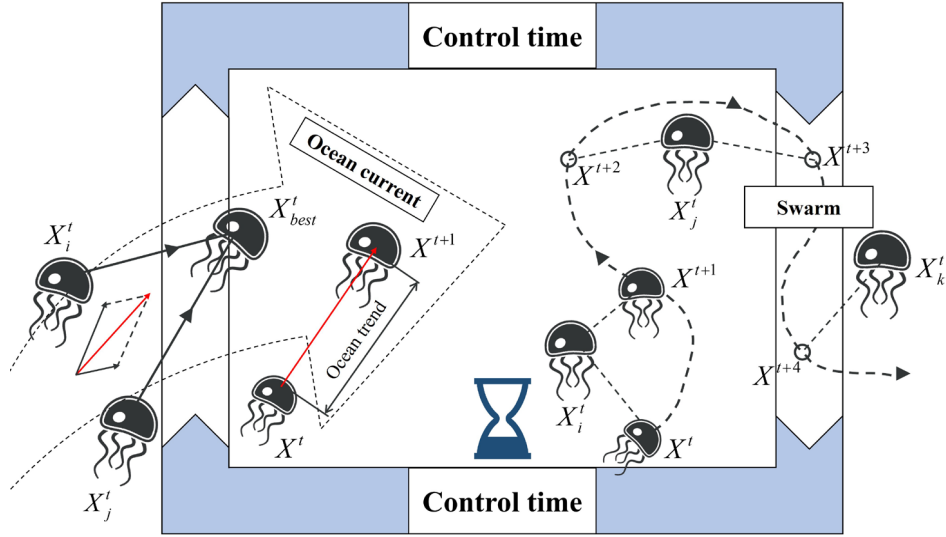

Fig. 5. Simulation of ocean current, swarm and time control

### 4.2.5 Levy Flight

In order to avoid the problem of the original JS falling into a local optimum and leading to premature convergence, the IJS introduces Levy flight. Levy flight is a strategy that enhances the global search capability. It can help the algorithm to jump out of the local optimum and explore a wider solution space [38]. New locations are updated based on Levy flight as follows:

$$Levy(\beta_{Levy}) = \frac{u}{|v|^{1/\beta_{Levy}}} \qquad (53)$$

$$X_i^{t+1} = X_i^t + \alpha \times Levy(\beta_{Levy}) \qquad (54)$$

where $\alpha$ is a scale parameter for the step size, which is a random number in the interval [0,1]; $u$ and $v$ are random numbers generated by normal distribution; $\beta_{Levy}$ is the stable distribution parameter of Levy flight, $\beta_{Levy} = 1.5$.

### 4.2.6 Differential Evolution

In this paper, the differential evolution (DE) strategy is introduced to enhance the global search capability and quality of solutions. Proposed by Storn and Price (1995) [39], DE is an optimization technique that mimics the process of genetic recombination in biological evolution. It generates new candidate solutions by combining the differences between the current optimal individual and other individuals. The specific steps are as follows:

**Step 1** (mutation): Two different individuals are randomly selected from the jellyfish population to generate a mutant vector. We use Equation (55) to realize the mutation:

$$V_i^t = X_i^t + \alpha^* \times (X_j^t - X_k^t) \tag{55}$$

where $X_j^t$ and $X_k^t$ are two different jellyfish individuals randomly selected from the population in the $t$th iteration; $\alpha^*$ is a scale parameter for the step size, which is a random number in the interval [0,1].

**Step 2** (crossover): We use a binary crossover strategy in this paper, where the crossover probability CR determines the extent to which the mutant vector is mixed with the current individual. The total dimension of each jellyfish individual is $D$, which is the sum of the dimensions of the bus scheduling scheme ($1 \times R^*$), the bus type allocation scheme ($1 \times K^*$), and the bus passenger capacity allocation scheme ($1 \times K^*$), i.e., $D = 1 \times (R^* + 2K^*)$. In the $t$th iteration, we introduce crossover between mutant vector and the current individual and form the trail vector $U_{i,d}^t$. For the $d$th dimension of $U_{i,d}^t$ ($d = 1, 2, ..., D$):

$$U_{i,d}^t = \begin{cases} V_{i,d}^t, & rand_{i,d} < CR \\ X_{i,d}^t, & rand_{i,d} \geq CR \end{cases} \tag{56}$$

where $V_{i,d}^t$ is the $d$th dimension of $V_i^t$; CR is set to 0.5 in this paper; $rand_{i,d}$ is a random number in the interval [0,1] to control the selection of the $d$th dimension of the $i$th individual in the crossover process.

**Step 3** (selection): Better individuals are selected for the next generation of the population by comparing the fitness of the trail vector with that of the current individual:

$$X_i^t = \begin{cases} U_i^t, & F(U_i^t) \geq F(X_i^t) \\ X_i^t, & F(U_i^t) < F(X_i^t) \end{cases} \tag{57}$$

### 4.2.7 Boundary conditions

When the jellyfish moves out of the bounded search area, it returns to the opposite bound based

on Equation (58):

$$X_{i,d}^{t+1} = \begin{cases} (X_{i,d}^t - U_{b,d}) + L_{b,d}, X_{i,d}^t > U_{b,d} \\ (X_{i,d}^t - L_{b,d}) + U_{b,d}, X_{i,d}^t < L_{b,d} \end{cases} \quad (58)$$

**4.2.8 Pseudocode of IJS**

The IJS, in addition to maintaining the advantages of the JS, can significantly improve the global search capability and convergence speed. We also show better performance and adaptability in the following experiments, with a total time complexity of $O(N \times (R^* + 2K^*) \times Max_{it})$. The pseudocode of the IJS for solving the proposed model is shown in Fig. 6.

```
Begin
    Define objective function F
    Set the search space, population size (N), and maximum iteration (Max_it)
    Initialize population of jellyfish X_i (i = 1, 2, ..., N) using tent chaotic map
    Calculate quantity of food at each X_i, F(X_i)
    Find the jellyfish at location currently with most food
    Initialize time: t=1
    Repeat
        For i=1: N do
            Calculate the time control C(t) using Eq. (52)
            If C(t)≥0.5: Jellyfish follows ocean current
                (1) New location of jellyfish is defined by Eq. (47)
            Else: Jellyfish moves inside a swarm
                If rand(0,1)>(1-C(t)): jellyfish exhibits passive motions
                    (1) New location of jellyfish is defined by Eq. (48)
                Else: Jellyfish exhibits active motions
                    (2) Determine direction of jellyfish using Eq. (49)
                    (3) New location of jellyfish is defined by Eq. (51)
                End if
            End if
            Levy Flight: New location of jellyfish is defined by Eq. (54)
            Differential Evolution: New location of jellyfish is defined by Eq. (57)
            Check boundary conditions and calculate quantity of food at new location
            Update the location of jellyfish
            Update the location of jellyfish currently with the most food
        End for i
        Update the time: t=t+1
    Until stop criterion is met (t > Max_it)
    Output the best results
End
```

Fig. 6. Pseudocode of the IJS

# 5. Numerical experiments

Computational experiments were conducted to (1) illustrate the good performance of the proposed algorithm for urban-rural bus networks, (2) compare the proposed PFSM with the separated passenger and freight transport, and (3) discuss the implications of the sensitivity of passenger-freight demand, maximum passenger travel time, minimum percentage of passenger capacity to vehicle capacity. The following experiment focuses on selecting evaluation indicators from three aspects, namely economy, efficiency and environment, to analyze the daily operational benefits of urban- rural bus systems. All algorithms were coded in MATLAB R2023a and executed on a Windows 10 computer equipped with AMD Ryzen 7 6800H, Radeon Graphics 3.20 GHz CPU, and 16 GB RAM.

## 5.1 Performance of IJS for urban-rural bus networks

### 5.1.1 Data description and parameter setting

In order to verify the effectiveness of the proposed IJS in solving optimization problems of urban-rural bus network under PFSM, genetic algorithm (GA), particle swarm optimization (PSO), grey wolf optimizer (GWO), and JS are selected as the benchmarks to compare their performances. In this section, the experiments were conducted using the simulated bus network consisting of seven bus lines proposed by Zeng et al [29], as shown in Fig. 7. Some lines have overlapping intervals. Both rural and urban bus terminals are adjacent to DCs, and parcels are loaded onto buses at stop 9 (rural DC) and unloaded at stop 10 (urban DC). The bus type for each line are classified as small, medium, and large. The relevant configuration parameters are shown in Table 1, and the departure timetable is shown in Table 2. We summarize the on-demand requests for passenger and cargo in Appendix Table A2.

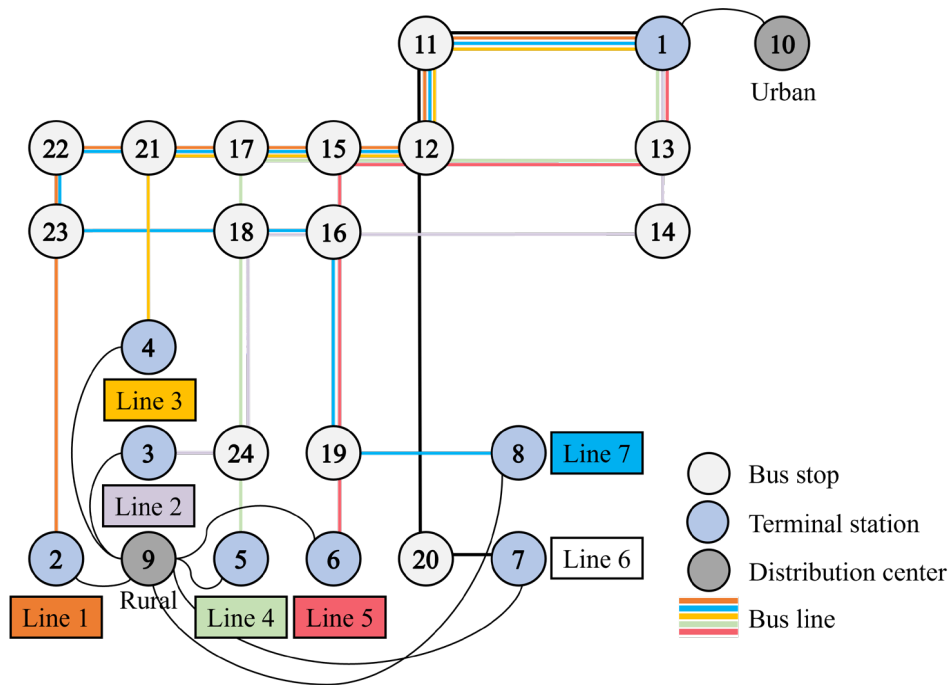

Fig. 7. Simulated urban-rural bus network

Table 1. Configuration parameters for different types of buses

| Vehicle type | Running cost $\phi^{km}$ (RMB/km) | Purchasing cost $\phi^{buy}$ (RMB/day) | Base fare of passenger service $\gamma_u$ (RMB/person) | Additional fare of passenger service $\eta^u$ (RMB/km) | Base fare of freight service $\gamma_f$ (RMB/item) | Additional fare of freight service $\eta_f$ (RMB/km) | Vehicle capacity $V$ (m³) |
|---|---|---|---|---|---|---|---|
| Small | 0.4 | 421.23 | | 0.17 | | | 13.56 |
| Medium | 0.6 | 506.23 | 2.5 | 0.22 | 1 | 0.05 | 20.3 |
| Large | 0.8 | 678.08 | | 0.27 | | | 30.6 |

Table 2 The departure timetable

| | Request ID | Departure | Arrival | Energy consumption |
|---|---|---|---|---|
| Line 1 (node 1 to node 2) | 1 | 15:00 | 16:50 | 30 kWh |
| | 2 | 17:00 | 18:50 | 30 kWh |
| Line 2 (node 1 to node 3) | 3 | 15:30 | 17:00 | 25 kWh |
| | 4 | 17:30 | 19:00 | 25 kWh |
| Line 3 (node 1 to node 4) | 5 | 16:30 | 18:00 | 25 kWh |
| Line 4 (node 1 to node 5) | 6 | 15:10 | 16:40 | 25 kWh |
| | 7 | 17:10 | 18:40 | 25 kWh |
| Line 5 (node 1 to node 6) | 8 | 16:50 | 18:20 | 25 kWh |
| Line 6 (node 1 to node 7) | 9 | 16:40 | 17:50 | 19 kWh |
| Line 7 (node 1 to node 8) | 10 | 17:20 | 20:00 | 44 kWh |

### 5.1.2 Benchmarks

We select GA, PSO, GWO, and JS as the benchmarks in the following experiments: (1) GA is a metaheuristic that mimics natural selection and the genetic principles (Holland, 1975) [40]; (2) PSO is a metaheuristic that simulates the behavior of animal cluster pattern (Kennedy & Eberhart, 1995) [41]; (3) GWO is a metaheuristic that simulates the hunting behavior of grey wolves (Mirjalili et al., 2014) [42]; (4) JS is metaheuristic that simulates the behavior of jellyfish in the ocean and is the original version of our proposed algorithm (Cho & Truong, 2021) [36].

### 5.1.3 Results and analysis

All algorithms have a population size of 100 and the maximum number of iterations of 150. We obtain the fitness and computational time for each algorithm for solving the optimal solution in the simulated network. Fig. 8 shows the iteration curves of the GA, PSO, GWO, JS and IJS.

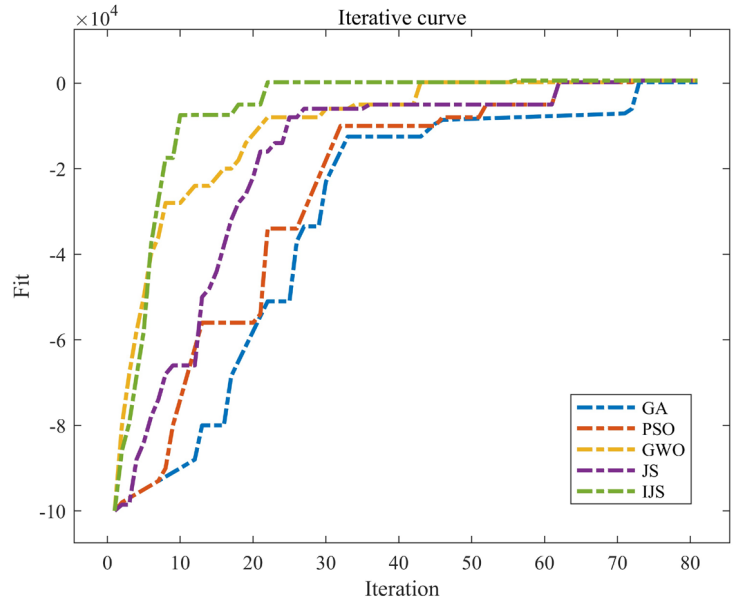

Fig. 8. Iteration curves of different algorithms

In the simulation experiments of the urban-rural bus network, IJS is the algorithm with the shortest solution time, the fastest convergence speed and the best solution quality among all algorithms. Specifically, except for the GA, PSO, GWO, JS, and IJS can find the same optimal solution by 80th iterations, which is better than the final solution obtained by the GA. Therefore, GA performs the worst both in terms of computational time and solution quality, which indicates that it is not suitable for dealing with the proposed problem. As an enhance version of the JS, the IJS is not only 51.98% shorter than the JS in terms of solution time and 3.66% shorter than the PSO, but also outperforms the other algorithms in terms of the quality of the solution at each iteration. It can be seen that the IJS has stronger global search capability and faster convergence speed than the JS, demonstrating excellent performance and adaptability in urban-rural transit network optimization problems. Therefore, the subsequent experiments use the IJS with optimal performance for model solving.

**5.2 Comparison between PFSM and separated passenger and parcel transport**

In order to quantify the advantages of PFSM in urban-rural bus operations, this paper compares the daily bus operation metrics of PFSM and separated passenger and parcel transport. In this paper, the separated passenger and parcel transport refers to medium-sized truck dispatching from the logistics system to transport parcels and medium-sized electric buses scheduling by the bus operators to serve passengers, both of which operate independently.

**5.2.1 Data description and parameter setting**

Based on the empirical study case of Zeng et al [29], the operation data of two urban-rural bus lines are collected in Yushe County, Shanxi Province, China, including information on the bus passenger-freight demand distribution of each line, timetable, vehicle scheduling, and bus fleet size. The two urban-rural bus lines depart from different rural starting terminals but share the same urban terminal and vehicle resources, as shown in Fig. 9. The rural terminals and the urban terminal are close to DC 11 and DC 10 at a distance of 3km and 2km, respectively. There is an ITSC in the middle section of both lines (ITSCs 4 and 16). Line 1 covers 8 stops, with a total length of 40.2km and an average single trip time of 70 minutes. Line 2 has a total length of 37.3km, 7 stops, and an average single trip time of 65 minutes. The distance matrix between stops is shown in Table 3. The freight lines are essentially the same as the bus lines, transporting parcels to and from rural DC 11 and urban DC 10, passing through ITSCs 4 and 16 on the way.

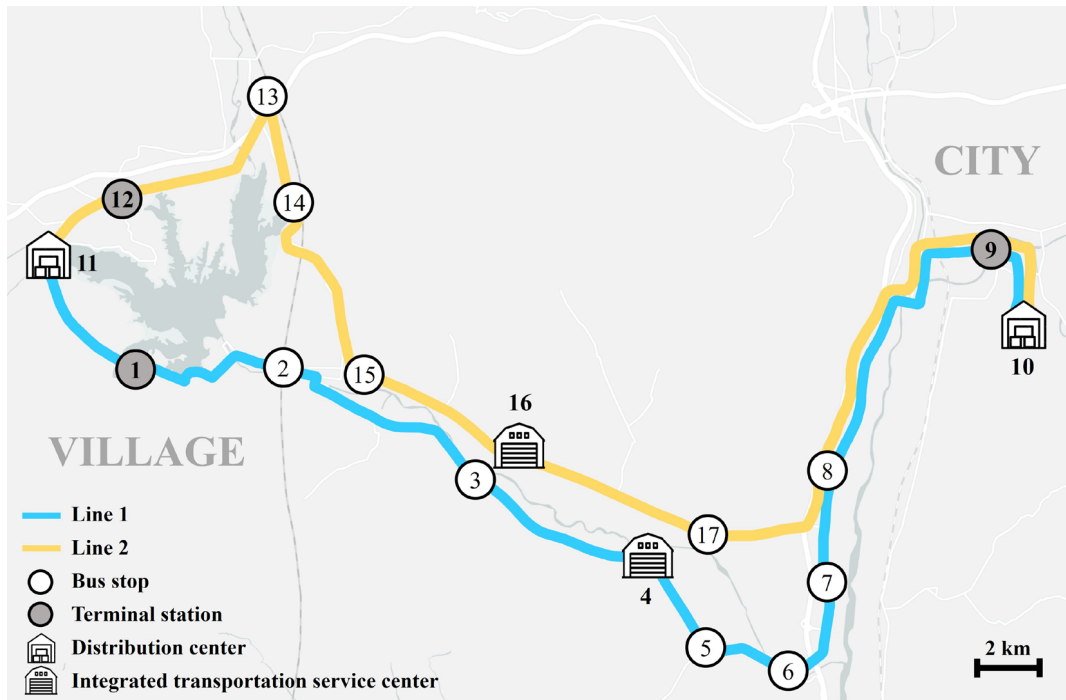

Fig. 9. Urban-rural bus lines in Yushe

Table 3 Bus line information of the urban-rural bus routes

| Bus line 1 | | | Bus line 2 | | |
| --- | --- | --- | --- | --- | --- |
| From | To | Distance (km) | From | To | Distance (km) |
| 10 | 1 | 3.00 | 10 | 12 | 3.00 |
| 1 | 2 | 6.12 | 12 | 13 | 3.63 |
| 2 | 3 | 6.58 | 13 | 14 | 3.19 |
| 3 | 4 | 4.90 | 14 | 15 | 4.40 |
| 4 | 5 | 2.16 | 15 | 16 | 3.49 |
| 5 | 6 | 3.60 | 16 | 17 | 5.83 |
| 6 | 7 | 2.60 | 17 | 8 | 4.51 |
| 7 | 8 | 2.00 | 8 | 9 | 12.28 |

| 8 | 9 | 12.28 | 9 | 11 | 2.00 |
| 9 | 11 | 2.00 | | | |

Due to the tidal characteristics of urban-rural commuting and the asymmetric distribution of passenger demand spatially and temporally, the following pattern exists in the timetable of the two bus lines: buses depart from the countryside to the city in the morning, from the city to the countryside in the afternoon, and stop in other off-peak hours, as shown in Table 4. The total daily demands for passengers and parcels are 300 and 1500, respectively. The OD demand information for passengers and parcels is shown in Appendix Table A3. Note that there are currently six buses serving two urban-rural bus lines, but the schedule does not fully utilize the fleet capacity because these buses also serve other bus lines.

Table 4 Bus line timetable

| Bus line 1 | | | Bus line 2 | | |
|---|---|---|---|---|---|
| Morning trip from node 1 to node 9 | | | Morning trip from node 12 to node 9 | | |
| No. | Departure | Arrival | No. | Departure | Arrival |
| 1 | 6:00 | 7:10 | 7 | 6:10 | 7:15 |
| 2 | 7:30 | 8:40 | 8 | 8:10 | 9:15 |
| 3 | 9:00 | 10:10 | 9 | 10:10 | 11:15 |
| Afternoon trip from node 9 to node 1 | | | Afternoon trip from node 9 to node 12 | | |
| 4 | 15:00 | 16:10 | 10 | 15:10 | 16:15 |
| 5 | 16:30 | 17:40 | 11 | 17:10 | 18:15 |
| 6 | 18:00 | 19:10 | 12 | 19:10 | 20:15 |

### 5.2.2 Results and analysis

For the separated passenger and parcel transport, both buses and trucks are medium-sized vehicles. According to the configuration parameter standard for different types of trucks in China, the vehicle capacity of a medium-sized truck is 17.28m$^3$, the cost of fuel consumption is RMB 1.6/km, the average purchasing cost is RMB 850/day, and the wage of a truck driver is RMB 500/day. The optimal scheme for vehicle allocation and scheduling under PFSM is shown in Table 5.

Table 5. Optimal scheme for vehicle allocation and scheduling under PFSM

| Bus number | Scheduling | Vehicle type | percentage of passenger capacity |
|---|---|---|---|
| 1 | No.1→ No.5 | Small | 100% |
| 2 | No.2→ No.4 | Small | 100% |
| 3 | No.3→ No.6 | Medium | 80% |
| 4 | No.7→ No.10 | Large | 50% |

| | | | | |
|---|---|---|---|---|
| 5 | No.8→ No.11 | Medium | | 90% |
| 6 | No.9→ No.12 | Large | | 70% |

As can be seen from Table 5, although the number of vehicles and the scheduling scheme remain unchanged, the vehicle type and the passenger capacity on each bus change significantly. As a result, the operational benefits change. A comparison of the daily operational metrics between the traditional separated passenger and parcel transport and PFSM is shown in Table 6.

Table 6. Comparison between PFSM and separated passenger and parcel transport

| Operational metric | Separated passenger and parcel transport | | PFSM | |
|---|---|---|---|---|
| | Truck | Bus | Truck | Bus |
| Fleer size | 4 | 6 | 0 | 6 |
| Total running distance (km) | 350.28 | 465.42 | 0 | 525.42 |
| $C_{idle}$ (RMB) | 55.8 | 26 | 0 | 69 |
| $C_{km}$ (RMB) | 589.48 | 390.6 | 0 | 389.21 |
| $C_{fix}$ (RMB) | 3400 | 3037.38 | 0 | 3211.08 |
| $C$ (RMB) | 4045.28 | 3453.98 | 3257.64 | 3669.29 |
| $E_u$ (RMB) | - | 991.82 | - | 991.82 |
| $Z$ (RMB) | - | **-2462.16** | - | **580.17** |
| $T_{cruise}$ (min) | - | 28.3 | - | 32.31 |
| $T_{det}$ (min) | - | 5.04 | - | 7.61 |
| $T_{dwell}$ (min) | - | 1.2 | - | 2.1 |
| $T_{wait}$ (min) | - | 4.15 | - | 4.2 |
| $\overline{T}$ (min) | - | **38.69** | - | **46.22** |

From Table 6, the traditional mode requires 4 medium-sized trucks and 6 medium-sized conventional buses, while PFSM requires only 6 buses to meet the passenger and freight transport requirements, which significantly reduces the number of vehicles required for freight services and effectively save the freight driver labor cost and fleet maintenance cost. PFSM has a differentiated vehicle allocation scheme, forming a combination of 'passenger-dominated' and 'freight-dominated' vehicles to adapt to the temporal and spatial mismatch of urban-rural demand. 4 small- and medium-sized buses are mainly responsible for passenger service, which are passenger-dominated. Two small buses serve passengers only, while the remaining two medium-sized buses share a small portion of the freight tasks. On the other hand, the two large buses make use of the more potential cargo space to undertake the main freight tasks of the original four medium-sized trucks, which are freight-dominated vehicles. Therefore, in this case with two DCs, PFSM can move the bus system into the black and make a profit of RMB 580.17 per day, which is a 328.45% increase in revenue compared with the traditional mode. It also reduces the total operating cost by 51.07%. Although the increase in freight volume affects passenger service, the average passenger travel time under PFSM only increases by 7.53 minutes, an increase of 19.46 %, which is within acceptable limits.

According to the Greenhouse Gas Protocol [43], the carbon emission factor of diesel fuel is 2.6765 kg $CO_2$/L. Since the fuel consumption of a medium-sized truck is 0.15 L/km, its $CO_2$ emissions are 0.40 kg/km. The daily freight transport mileage of the conventional mode is 350.28 km, which results in 140.60 kg of $CO_2$ emissions. Although electric buses use electricity and do not directly generate carbon emissions, electricity is mainly generated by coal combustion, which generates carbon emissions of 0.7967kg/kWh [44]. In terms of energy saving and emission reduction, the PFSM replaces medium-sized trucks to undertake freight services by optimizing the configuration of the electric bus type and capacity allocation. The daily electricity consumption is 394.06kWh and the $CO_2$ emissions are 313.95kg, which is 25.02% lower than the separated passenger and parcel transport. In summary, in this case with two DCs, PFSM can transport about 547,500 express deliveries per year, saving 19,178 L of freight fuel consumption and reducing $CO_2$ emissions by 38.23 tons. If this mode is extended to the whole China, it can save about 3,317.79 million liters of fuel consumption and directly reduce $CO_2$ emissions by about 6,614,200 tons, which is equivalent to reducing the annual $CO_2$ emissions from motor vehicles of China in 2023 by 1.06%. This is a significant contribution to reducing carbon emissions and improving air quality, providing an effective solution to green transport.

**5.3 Sensitivity Analysis**

To demonstrate the reliability and validity of optimal solutions with respect to changes in the values of parameters, this section further investigates the sensitivity of passenger-freight demand, maximum passenger travel time, and minimum percentage of passenger capacity to vehicle capacity based on the case in Section 5.2.

**5.3.1 Passenger-freight demand**

Passenger-freight demand is an important factor that affects the operating profit and passenger travel time. Fig. 8 displays the effect of different passenger and freight demands on the operating profit and passenger travel time. The passenger demand takes the values of {100, 300, 500, 700, 900} and the freight demand takes the values of {500, 1000, 1500, 2000, 2500}, representing the change process of passenger and freight demand levels from low to high. The OD demand values between stops are modified based on the ratio of the above total demand to the original total demand in Appendix Table A2.

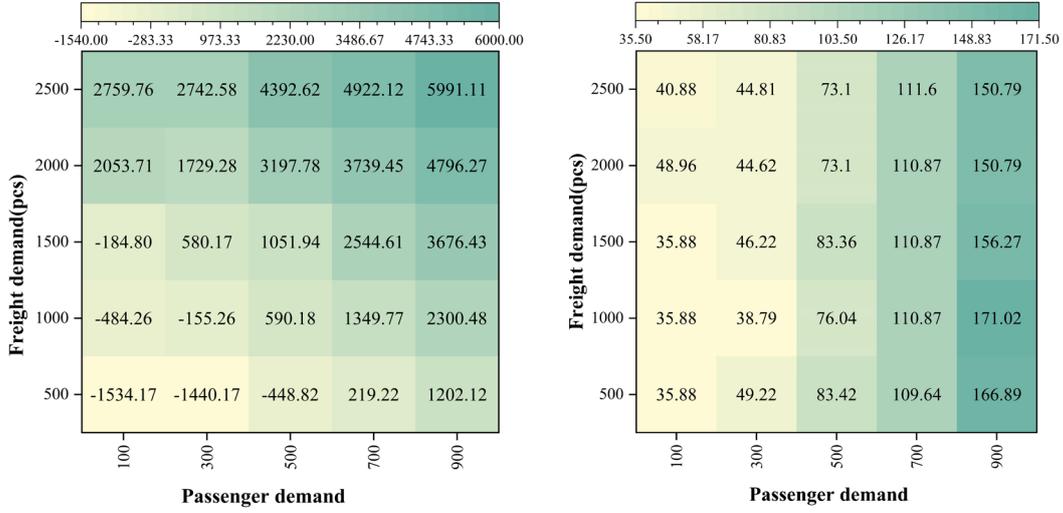

(a) operating profit (RMB/day)　　(b) passenger average travel time (min)

Fig. 10. Operating profit and passenger average travel time under PFSM with different passenger-freight demands

From Fig. 10 (a), for the operating profit $Z$, it can be seen that bus operations are prone to losses when the passenger-freight demand is low, indicating that PFSM may not be suitable for scenarios with low demand. In particular, the minimum value of -1534.17 RMB/day is reached when the passenger-freight demand is {100, 500}. As the passenger-freight demand increases to {100, 2000}, {300, 1500}, {500, 1000}, or {700, 1000}, the bus system is moved into the black. The above passenger-freight demands may be the minimum criteria to promote the profitability of the bus system under PFSM. When the passenger-freight demand is {900, 2500}, it reaches the maximum value of 5991.1 RMB/day. The color change of the heat map clearly shows the gradual increase of profit from negative (loss) in the lower left corner to positive (profit) in the upper right corner, which indicates that the total profit increases with the increase of passenger-freight demand. Although total profit increases whenever either passenger or freight demand increases, there is a difference in the extent to which passenger and freight demands contribute to operating profit.

From Fig. 10 (b), for the passenger average travel time $\bar{T}$, it can be seen that the minimum value of 35.88 min is reached when the passenger-freight demand is {100, 500}, {100, 1000}, and {100, 1500}, respectively. The maximum value of 171.02 min is reached when the passenger-freight demand is {900, 1000}, resulting in a large number of secondary passenger detentions. Different from $Z$, the heat map has a more pronounced color change horizontally than vertically. The similar colors in the same column suggest that passenger demand has a greater impact on $\bar{T}$ than freight demand, showing a nonlinear growth. Under the same passenger demand scenario, the range of $\bar{T}$ values is relatively stable compared with Fig. 10 (a). This indicates that an increase in the number of passengers significantly increases the boarding/alighting time and the passenger detention time, thus seriously affecting $\bar{T}$. This effect tends to intensify with the incremental increase in passenger demand and the efficiency of PFSM will face more challenges. Besides, in some combinations of passenger-freight demand, an increase in freight demand reduces $\bar{T}$ instead. This is because under higher freight demand, fleet configuration tends to select medium- and large-sized buses in order to satisfy both passenger and freight demands, thus expanding the space for carrying more passengers and parcels, and improving the efficiency and the passenger travel experience.

In order to further explore the profit contribution differences between passengers and cargoes under different passenger demand levels, we introduce profit contribution ratio (PCR), space-based profit contribution ratio (SPCR), and income equivalence ratio (IER) as the evaluation indicators. PCR denotes the number of parcels transported whose profit is equal to the profit made by transporting one passenger; IER indicates the ratio of the number of passengers to the number of parcels to be transported to obtain the same revenue; and SPCR represents the profit ratio of using the same capacity to serve passengers and freight separately. In general, the capacity of a seat allows for stacking of 10 small parcels.

Table 7. Difference in profit contribution between passengers and parcels at different passenger demand levels

| Passenger demand level | PCR | IER | SPCR |
|---|---|---|---|
| Low (100) | 1:0.59 | 1:1.693 | 1:5.91 |
| Medium (500) | 1:0.50 | 1:2.0107 | 1:4.98 |
| High (900) | 1:0.46 | 1:2.1988 | 1:4.57 |

From the above table, the freight service profit is greater than the passenger service profit in any passenger demand scenario, and the gap widens gradually with the growth of passenger demand. Based on PCR and IER, the unit profit contribution of freight transport is 0.46-0.59 times that of passenger transport. In terms of SPCR, the economic value of transporting parcels in the same space is about 4.57-5.91 times that of passenger transport. In summary, the economic advantage of freight service is more significant under PFSM, and that advantage is greater as passenger demand grows, probably because of the saturation of vehicle occupancy in the high passenger demand scenario and the need to increase vehicle capacity to serve both passengers and parcels.

**5.3.2 Maximum passenger travel time**

The maximum passenger travel time $T_{max}$ is an important factor to ensure the operational efficiency of urban-rural bus routes, and it plays a role in regulating the conflict between economic efficiency and service level of the bus operation scheme. If the value of $T_{max}$ is too large, the scheme may neglect passengers travel experience and satisfaction, thus reducing the attractiveness and operational efficiency of urban-rural buses; if the value is too small, it may limit the optimization space of the bus operation scheme and reduce the economic and environmental benefits of bus operation, even causing the fail to obtain a feasible solution [23]. Fig. 11 summarizes the economic and time metrics for the operation of urban-rural bus routes under different values of $T_{max}$.

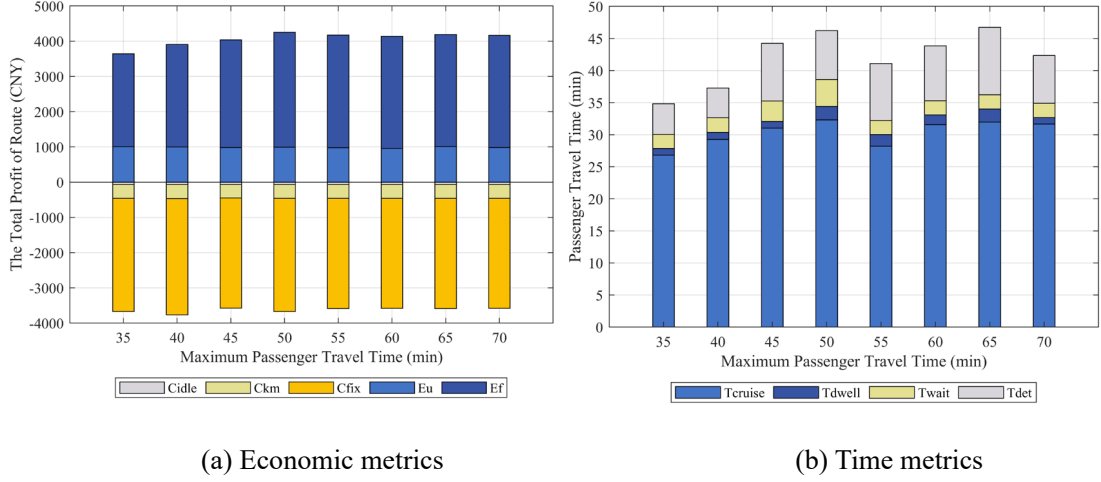

(a) Economic metrics        (b) Time metrics

Fig. 11. Economic and time metrics for the operation of urban-rural bus routes under different values of $T_{max}$

It can be seen that the total operating profit $Z$ increases with the increase of $T_{max}$, and tends to stabilize when $T_{max}$ is greater than 50 minutes. When $T_{max}$ = 35, $Z$ = -26.53 RMB/day, indicating that the bus system incurs a slight loss in order to fully guarantee the efficiency of passenger travel. This is mainly because the profit of transporting parcels with limited vehicle capacity is higher than the profit of carrying passengers under PFSM, but in order to satisfy the $T_{max}$ constraint, the optimization model guarantees the adequate passenger supply by reducing the freight capacity. The reduction of the revenue from freight services makes a slight loss. The profit increases rapidly with the gradual increase of $T_{max}$. When $T_{max} \geq 50$, the value of $Z$ fluctuates within the range of [553.62, 605.41], which is 122%-124% more than the profit of the traditional separated passenger and parcel transport. This implies that the operating profit of the bus system has reached a nearly saturated state, where the capacity allocation of passenger and freight services is reasonable to meet the corresponding demand.

In terms of average passenger travel time $\bar{T}$, when $T_{max}$ is less than 40 minutes, $\bar{T}$ approaches $T_{max}$. As $T_{max}$ increases, $\bar{T}$ is no longer close to $T_{max}$ and stabilizes at 41-46 minutes, an increase of 4%-20% compared to the conventional mode. This shows that the optimal scheme has balanced the benefits between passengers and the operator within the acceptable range of $T_{max}$, rather than sacrificing the interests of passengers endlessly. In order to balance the economic benefits and the level of service, it is recommended that $T_{max}$ is set to 50 min.

### 5.3.3 Minimum percentage of passenger capacity to vehicle capacity

For the minimum percentage of passenger capacity to vehicle capacity on a bus ($\lambda_{min}$), choosing the appropriate value allows the optimal scheme to balance operational benefits and service quality. If $\lambda_{min}$ is too small, it may lead to an oversupply of passengers and reduce the level of bus service; if $\lambda_{min}$ is too large, it may affect the freight demand and reduce the revenue. Fig. 12 shows the economic and time metrics of urban-rural bus operations under different values of $\lambda_{min}$.

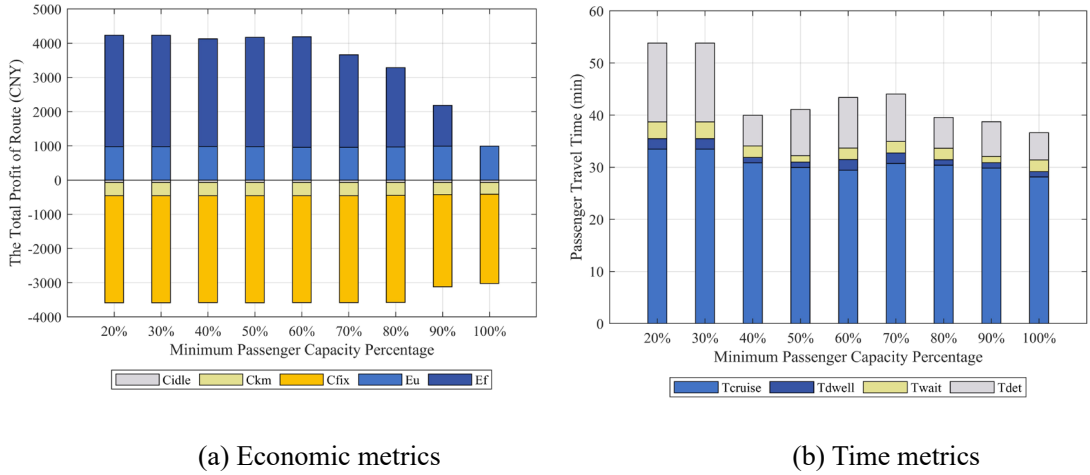

(a) Economic metrics    (b) Time metrics

Fig. 12. Economic and time metrics for the operation of urban-rural bus routes under different values of $\lambda_{min}$

As $\lambda_{min}$ increases, the total operating profit $Z$ is stable at around 4200 RMB/day and then decreases, while $\overline{T}$ has a fluctuating trend of decreasing, then increasing and decreasing again. In terms of $Z$, when $\lambda_{min} < 40\%$, $Z$ reaches the maximum value of 653.17 RMB/day; when 40% $\leq \lambda_{min} \leq 60\%$, $Z$ slightly decreases but is stable in the interval [543.71, 607.01], increasing by 122%-125% than the profit of separated passenger and parcel transport. When $\lambda_{min} > 60\%$, $Z$ decreases rapidly as $\lambda_{min}$ increases and begins to turn from profit to loss. The main reason is that buses with a small $\lambda_{min}$ can satisfy most freight demands with considerable economic benefits. As $\lambda_{min}$ increases, the remaining freight capacity on buses shrinks, resulting in a decline in revenue from freight services.

In terms of average passenger travel time $\bar{T}$, when $\lambda_{\min}$ increases from 20% to 30%, $\bar{T}$ remains unchanged, mainly because the freight supply still exceeds demand in this condition and the reduction in freight capacity does not change the total dwell time. When $\lambda_{\min}$ grows from 30% to 40%, the equilibrium between freight supply and demand is disrupted, and the compression of freight capacity leads to a reduction in the number of parcels carried, thus reducing the delay caused by freight service, with a significant decrease of $\bar{T}$. However, when 40% ≤ $\lambda_{\min}$ ≤ 70%, an increase in $\lambda_{\min}$ increases $\bar{T}$ by 3%-12% compared with the traditional mode. This is mainly due to the fact that the increase in the number of passengers carried leads to the additional dwell time while the expanding passenger capacity accommodates more passengers. When $\lambda_{\min}$ >70, $\bar{T}$ decreases to a minimum value of 36.64 min as $\lambda_{\min}$ increases. One reason is that the passenger supply is greater than the demand and the expansion of the passenger capacity no longer increases the number of passengers. The other reason is that the reduction of the freight capacity results in the decreasing number of parcels, thus reducing the travel delay by freight services. In summary, we recommend $\lambda_{\min} \in [40\%, 60\%]$ to ensure the efficiency of passenger travel at the simultaneous maximization of operating profits.

**5.4 Policy implications and applications**

The policy values of this study are embodied in the collaboration to promote economies in transition, low-carbon development, and rural revitalization, with the following connotations:

(i) Economic dimension: constructing a virtuous cycle of two-way feedback between urban-rural logistics cost control and revenue generation from public transport

The quantitative analysis on the case of Shanxi shows that this study can reduce the cost of rural end-to-end logistics by 19.47%. The paths to realizing its economic value consist of 1) the effect of economies of scale by taking advantage of bus network coverage to consolidate fragmented freight demand and reduce rural logistics costs as a lower share of GDP; 2) the revenue-sharing mechanism by setting appropriate pricing strategies to achieve synergies between passenger and freight services and increase the operating benefits [45]; 3) the resource reconfiguration by forming a combination of passenger-dominated and freight-dominated buses and reducing the logistics fleet size from rural areas to the urban core, in order to adapt to the temporal and spatial mismatch of urban-rural demand, increase vehicle turnover, and reduce the freight cost. Therefore, we suggest establishing a system of 'public transport for freight transport franchising with revenue feedback', whereby part of the logistics profits will be directed to the improvement of bus service quality such as increasing the bus frequency in peak hours and optimizing the bus stop layout, forming a sustainable urban-rural transport economic ecosystem.

(ii) Environmental dimension: pioneering pathways for optimizing low carbon transport systems

PFSM has significant advantages of carbon reduction over traditional transport systems, including: 1) at the level of energy structure optimization, the replacement of diesel trucks by electric buses reduces 19.12 tons of $CO_2$ emissions per year on a single route, which is equivalent to the carbon sink capacity of planting 875 trees; 2) at the level of energy efficiency improvement in transport, the flexible allocation strategy for freight capacity improves the space utilization of vehicle capacity by 34.8%, thereby reducing the energy consumption per unit transport turnover of freight. Therefore, we can refer to the carbon accounting system for transportation, promoting the carbon credit trading mechanism based on the life cycle assessment (LCA), and provide financial subsidies for carbon reduction to the public transport enterprises implementing PFSM to accelerate the achievement of carbon peaking and carbon neutrality goals in the transport sector.

(iii) Social dimension: shaping a new framework of integrated urban-rural development

PFSM has significant values of social inclusion: 1) in terms of service accessibility, it shortens the average transport time for urban-rural express services and increases e-commerce penetration in remote rural areas; 2) in terms of job creation, it not only avoids massive layoffs due to a reduction in the scale of bus operations, but also creates more employment opportunities in rural areas such as cargo handling and management, thus contributing to an increase in the average annual per capita income of rural residents; 3) in terms of travel security, the widespread application of PFSM helps public transport companies to develop more bus routes and provide better travel services, thereby benefiting rural residents who face some travel problems including the single travel mode and weak transport accessibility in rural areas [47]. Therefore, we can construct an investment system including basic subsidies from the central government, special bonds support from local governments, and the participation by social capital to promote the high-quality development of urban-rural economic integration.

## 6. Conclusion

To further expand the optimization theory of urban-rural bus operation under PFSM and maximize its economic and environmental benefits, this paper proposes an 'economy-efficiency-low-carbon' -oriented resource reconfiguration strategy. Considering the collaborative optimization of multi-type bus scheduling and dynamic vehicle capacity allocation for urban-rural bus routes, we formulate a bilevel optimization model to investigate this strategy, and then a novel improved jellyfish search algorithm is designed to solve the transforming model by the entropy weight method. Two case studies of a simulated bus network and two urban-rural bus lines in Shanxi Province, China, are used to examine the effectiveness of the proposed methodology in improving the economic benefits of urban-rural bus systems, reducing carbon emissions, and guaranteeing the level of bus service, respectively. We also obtain policy insights and give recommendations to the government on how to regulate the urban-rural bus operation under PFSM.

The proposed model can be further extended in future research. For example, the proposed model is based on the known OD demand, but there exist large uncertainties in passenger-freight transport demand in the actual operation, which has a certain impact on the reliability of the model application. Besides, due to space limitation, we do not include the timetable optimization in the model. However, there is an interactive feedback relationship between timetable setting and vehicle scheduling, which will significantly raise the requirements for the complexity and accuracy of the model solution. Therefore, our future research work will be devoted to investigating the influence of passenger-freight demand fluctuations on the travel delays and optimize bus timetable and vehicle scheduling cooperatively to enhance the theoretical and practical values.

## Data availability

The datasets generated and analysed during the current study will be made publicly available via a persistent URL upon formal publication of the article, in accordance with funder and institutional requirements for data sharing.

## Credit authorship contribution statement

**Jiabin Wu:** Conceptualization, Methodology, Funding acquisition, Writing-Original draft. **Zijian Huang:** Data curation, Validation, Visualization. **Linhong Wang:** Supervision, Writing-Reviewing and Editing, Methodology. **Yiming Bie:** Supervision, Funding acquisition, Writing-Reviewing and Editing. **Yuting Ji:** Formal analysis, Software. **Jun Gong:** Investigation, Data processing.

## Declaration of Competing Interest

The authors declare that they have no known competing financial interests or personal relationships that could have appeared to influence the work reported in this paper.

## Acknowledgments

This research is supported by the Humanities and Social Sciences Youth Foundation, Ministry of Education [grant number 23YJC630189]; Guangdong Basic and Applied Basic Research Foundation [grant number 2022A1515111039]; National Natural Science Foundation of China [grant number 52472359 and 72471102]; Transportation Innovation and Development Support Project of Jilin Province [grant number 2023-1-13]; Jilin University Lixin Excellent Young Teacher Training Program [grant number 4190805LX104].

of transportation electrification: Urban buses. Energy policy, 148, 111921.

[46] https://sustainable-bus.com/electric-bus/electric-bus-public-transport-main-fleets-projects-around-world/ (2021)

[47] Xue Bing, Zhi-Chun Li, Xiaowen Fu, Optimization of urban-rural bus services with shared passenger-freight transport: Formulation and a case study, Transportation Research Part A: Policy and Practice, 2025, 192, 104355. https://doi.org/10.1016/j.tra.2024.104355.

# Appendix

Table A1. List of notations

| Abbreviations | Definition | | |
|---|---|---|---|
| | | $x_k^r$ | Binary variable; equals 1 if bus $k$ serves run $r$, and 0 otherwise |
| PFSM | Passenger-Freight Shared Mobility | $y_k^n$ | Binary variable; equals 1 if the type of bus $k$ is type $n$, and 0 otherwise |
| DC | Distribution centers | $F$ | System benefit |
| ITSC | Integrated transportation service centers | $Z$ | Total operating profit |
| IJS | Improved Jellyfish Search Algorithm | $C_{km}$ | Total running cost |
| JS | Jellyfish Search Algorithm | $C_{dwell}$ | Total dwelling cost |
| GA | Genetic Algorithm | $C_{fix}$ | Total purchasing cost |
| PSO | Particle Swarm Optimization | $C_{toll}$ | Total toll cost |
| GWO | Grey Wolf Optimizer | $E_u$ | Transit passenger fare revenue |
| **Parameters** | **Definition** | $E_f$ | Revenue of freight service |
| $\phi_n^{km}$ | Running cost per kilometer of bus type $n$ | $T$ | Total passenger travel time |
| $\phi_{dwell}$ | Value of dwell time | $T_{cruise}$ | Total in-vehicle travel time |
| $\phi_n^{buy}$ | Purchasing cost of bus type $n$ | $T_{dwell}$ | Total dwell time |
| $T_{max}$ | Maximum passenger travel time | $T_{wait}$ | Total passenger waiting time |
| $\eta_n^u$ | Charged thereafter per kilometer of bus type $n$ for passenger service | $T_{det}$ | Total passenger detention time |
| $\eta_f$ | Charged thereafter per kilometer for freight service | $T_{i,j}^r$ | In-vehicle travel time between $(i, j)$ of the bus serving the $r$th run |
| $v_f$ | Free-flow travel time of roadway segments | $\varepsilon_{i,j}^r$ | Stochastic perturbation term under the normal distribution $N(0, \sigma^2)$ |
| $L_u$ | Initial distance included in the base fare of passenger service | $P_{r,i}^{det}$ | Number of stranded passengers after the $r$th run leaves stop $i$ |
| $L_f$ | Initial distance included in the base fare of freight service | $P_{r,i}^{wait}$ | Number of passengers waiting for the $r$th run at stop $i$ |
| $\gamma_u$, $\gamma_f$ | Base fares of passenger and freight services | $P_{r,i}^{re}$ | Number of seats remaining when the bus serving the $r$th run arrives at stop $i$ |
| $L_r$ | Running distance of the $r$th run in km | $o_{r,i}^u$, $d_{r,i}^u$ | Number of passengers boarding and alighting at stop $i$ on the bus serving the $r$th run |

| Symbol | Description | Symbol | Description |
|---|---|---|---|
| $V_n$ | Vehicle capacity of bus type $n$ in m³ | $o_{r,i}^{f}$, $d_{r,i}^{f}$ | Number of parcels loading and unloading at stop $i$ on the bus serving the $r$th run |
| $r_{i,j}$, $p_{i,j}$ | Performance value and weight of the $i$th alternative based on $j$th evaluation index | $d_{i,j}^{r}$, $q_{i,j}^{r}$ | Number of passengers and the amount of cargo transported from stop $i$ to stop $j$ in the $r$th run |
| $O$ | Number of alternatives | $\Delta_{r,m}^{u(i,j)}$ | Dwell times of a bus serving the $r$th run at stop $m$ |
| $K^*$ | Total fleet size | $\Delta_p$ | Average delay of a boarding/alighting passenger in seconds |
| $Max_{it}$ | Maximum number of iterations for the IJS | $\Delta_f$ | Average delay of a loading/unloading parcel in seconds |
| $U_b$, $L_b$ | Upper bound and lower bound of search spaces | $\Delta_{toll}$ | Road toll of a bus |
| $\beta_d$ | Distribution coefficient | $L_{toll}$ | Toll mileage of a bus route in km |
| $\gamma'$ | Motion coefficient related to the length of motion around the locations of the jellyfish | $a_{i,j}^{r}$, $Q_{i,j}^{r}$, $C_{i,j}^{r}$ | Free-flow travel time, demand volume, and capacity of the bus serving the $r$th run between $(i,j)$ |
| $C_0$ | Threshold constant to regulate passive and active movements of the jellyfish | $\hat{T}_{i,j}^{r}$, $\mu_{i,j}^{r}$, $\sigma_{i,j}^{r}$, $\Omega(\cdot)$ | Time budget value, mean, standard deviation, and distributing function of in-vehicle travel time of the bus serving the $r$th run between $(i,j)$ |
| $\alpha$, $u$, $v$ | Parameters for Levy flight | $R_{i,j}^{r}$ | Probability that the actual in-vehicle travel time from stops $i$ to $j$ of the bus serving the $r$th run is not greater than the travel time budget |
| $\alpha^*$ | Scale parameter for the step size in the mutation step of the differential evolution | $t_{r,i}$, $t'_{r,i}$ | Times when the $r$th run arrives at stop $i$ and the $(r-1)$-th run departs from stop $i$ |
| $R^*$ | Daily number of bus runs | $X_i$ | Initial location of the $i$th jellyfish |
| $D$ | The total dimension of a scheme for bus allocation and scheduling | $X_i^t$, $X_{i,d}^t$ | Location of the $i$th jellyfish in the $t$th iteration and the $d$th dimension of $X_i^t$ |
| $\mu_{tent}$ | Parameter that controls the distribution pattern of the tent chaotic sequence | $U_{i,d}^t$ | $d$th dimension of trail vector of the $i$th individual in the $t$th iteration |
| $\beta_{Levy}$ | Stable distribution parameter of Levy flight | $X_{best}^t$ | Best location in the current jellyfish population |
| $CR$ | Crossover probability of the differential evolution | $\bar{\mu}$ | Mean location of all jellyfish in the population |
| $rand_{i,d}$ | Random number to control the selection of the $d$th dimension of the | $V_i^t$, $V_{i,d}^t$ | Mutant vector of the $i$th individual in the $t$th iteration and the $d$th dimension of $V_i^t$ |
43

| Symbol | Definition | Symbol | Definition |
|---|---|---|---|
| | $i$th individual in the crossover process | | |
| $\gamma$ | Confidence level of in-vehicle travel time | $u(i,j)$ | Equals 1 if stop $i$ locates upstream of stop $j$ in the stop sequence; and -1 otherwise |
| $n^*$ | Number of evaluation indices | $\tilde{S}$, $\tilde{K}$ | Skewness and kurtosis of the in-vehicle travel time |
| $\alpha_i$ | Toll rate for class $i$ | $w_j$ | Weighting factor of the $j$th index |
| $\beta_k^n$ | Seat number of bus $k$ with type $n$ | $E_j$ | Information entropy of the $j$th evaluation index |
| $\bar{v}_p$ | Average spatial volume of a seat in a bus in m³ | **Auxiliary Symbols** | **Definition** |
| $S$ | Set of bus stops | $C(t)$ | Time control function $C(t)$ |
| $\beta$, $z$ | Calibration parameters of the improved BPR function | $\lfloor \cdot \rfloor$ | Floor function |
| $R$, $N$, $K$ | Sets of bus runs, vehicle types, and bus numbers | $\Phi_\gamma$ | $\gamma$-quantile of a standard normal distribution |
| $\beta_j$ | Seat number threshold for different classes | $D_i(t)$ | Cumulative passenger arrival function at stop $i$ |
| **Variables** | **Definition** | | |
| $\lambda_k$ | Percentage of passenger capacity to vehicle capacity on bus $k$ | | |
| $\delta_i$ | Binary variable; equals 1 if stop $i$ allows for cargo loading/unloading, and 0 otherwise | | |

Table A2 On-demand request for passenger (P) and cargo (C)

| No. | From | To | Pickup time window | Load (P/C) | No. | From | To | Pickup time window | Load (P/C) |
|---|---|---|---|---|---|---|---|---|---|
| 1 | 1 | 2 | [14:55, 15:25] | 30 P | 21 | 12 | 24 | [16:40, 17:50] | 30 P |
| 2 | 11 | 22 | [15:10, 15:40] | 30 P | 22 | 1 | 6 | [16:50, 17:20] | 30 P |
| 3 | 21 | 23 | [16:00, 16:30] | 30 P | 23 | 12 | 19 | [17:10, 17:40] | 30 P |
| 4 | 1 | 23 | [16:50, 15:20] | 30 P | 24 | 15 | 6 | [17:10, 17:40] | 30 P |
| 5 | 15 | 2 | [17:20, 17:50] | 30 P | 25 | 1 | 7 | [16:40, 17:50] | 30 P |
| 6 | 17 | 22 | [17:35, 18:05] | 30 P | 26 | 1 | 12 | [16:40, 17:50] | 30 P |
| 7 | 1 | 3 | [15:25, 15:55] | 30 P | 27 | 11 | 20 | [17:00, 17:30] | 30 P |
| 8 | 13 | 16 | [15:35, 16:00] | 30 P | 28 | 1 | 18 | [17:20, 18:00] | 30 P |
| 9 | 14 | 24 | [15:35, 16:00] | 30 P | 29 | 17 | 19 | [18:00, 18:30] | 30 P |
| 10 | 1 | 24 | [17:30, 18:00] | 30 P | 30 | 18 | 8 | [18:00, 18:30] | 30 P |
| 11 | 18 | 3 | [18:20, 18:50] | 30 P | 31 | 10 | 9 | [15:00, 15:30] | 250 C |
| 12 | 16 | 24 | [17:55, 18:25] | 30 P | 32 | 10 | 9 | [15:10, 15:40] | 250 C |

| 13 | 1 | 4 | [16:25, 16:55] | 30 P | 33 | 10 | 9 | [15:30, 16:00] | 250 C |
| 14 | 12 | 4 | [17:00, 17:30] | 30 P | 34 | 10 | 9 | [16:30, 17:00] | 250 C |
| 15 | 15 | 21 | [17:25, 17:55] | 30 P | 35 | 10 | 9 | [16:40, 17:10] | 250 C |
| 16 | 1 | 5 | [15:10, 15:40] | 30 P | 36 | 10 | 9 | [16:50, 17:20] | 250 C |
| 17 | 1 | 15 | [15:10, 15:40] | 30 P | 37 | 10 | 9 | [17:00, 17:30] | 250 C |
| 18 | 12 | 18 | [15:20, 15:50] | 30 P | 38 | 10 | 9 | [17:10, 17:40] | 250 C |
| 19 | 1 | 24 | [17:10, 17:40] | 30 P | 39 | 10 | 9 | [17:20, 17:50] | 250 C |
| 20 | 17 | 5 | [17:40, 18:10] | 30 P | 40 | 10 | 9 | [17:30, 18:00] | 250 C |

Table A3 On-demand request for passenger (P) and cargo (C)

| No. | From | To | Pickup time window | Load (P/C) |
|---|---|---|---|---|
| 1 | 3 | 9 | [6:10, 8:10] | 15 P |
| 2 | 12 | 9 | [6:30, 8:30] | 15 P |
| 3 | 16 | 8 | [7:00, 8:00] | 15 P |
| 4 | 8 | 9 | [7:25, 9.25] | 15 P |
| 5 | 10 | 4 | [7:45, 9:45] | 15 P |
| 6 | 4 | 8 | [8:00, 10:00] | 15 P |
| 7 | 15 | 9 | [8:35, 11:35] | 15 P |
| 8 | 8 | 9 | [9:10, 11:10] | 15 P |
| 9 | 2 | 8 | [9:20, 11:20] | 15 P |
| 10 | 14 | 17 | [10:10, 12:10] | 15 P |
| 11 | 9 | 12 | [15:10, 16:40] | 15 P |
| 12 | 4 | 10 | [15:20, 16:20] | 15 P |
| 13 | 9 | 15 | [16:00, 17:10] | 15 P |
| 14 | 9 | 3 | [16:30, 18:00] | 15 P |
| 15 | 17 | 12 | [16:40, 17:50] | 15 P |
| 16 | 8 | 4 | [16:45, 17:45] | 15 P |
| 17 | 9 | 1 | [17:00, 18:30] | 15 P |
| 18 | 8 | 16 | [17:40, 19:50] | 15 P |
| 19 | 9 | 8 | [18:00, 19:00] | 15 P |
| 20 | 8 | 1 | [18:30, 19:30] | 15 P |
| 21 | 10 | 11 | [6:00, 8:00] | 150 C |
| 22 | 10 | 11 | [6:10, 8:10] | 120 C |
| 23 | 10 | 4 | [7:30, 9:30] | 120 C |
| 24 | 4 | 11 | [8:10, 10:10] | 120 C |
| 25 | 10 | 16 | [9:00, 11:00] | 120 C |
| 26 | 16 | 11 | [10:10, 12:10] | 120 C |
| 27 | 11 | 10 | [15:00, 17:30] | 150 C |
| 28 | 11 | 10 | [15:10, 17:30] | 120 C |
| 29 | 11 | 4 | [16:30, 18:30] | 120 C |
| 30 | 4 | 10 | [17:00, 18:30] | 120 C |

| 31 | 11 | 16 | [18:00, 19:30] | 120 C |
| 32 | 16 | 10 | [19:00, 19:30] | 120 C |